\documentclass[final]{elsarticle}
\usepackage{lscape}
\usepackage[english]{babel}
\usepackage{url}
\usepackage{amsmath}
\usepackage{a4}
\usepackage{amssymb, amsmath}           % Mathematical packages
\usepackage{amsthm}
\usepackage{mathrsfs}
\usepackage{epsfig}
\usepackage{color,graphicx}             % eps compatibility
\usepackage{setspace}
\usepackage{caption}
\usepackage{algorithm} 
\usepackage[noend]{algpseudocode}
\graphicspath{{fig/}}

\newtheorem{assumption}{Assumption}

\newtheorem{remark}{Remark}
\usepackage{subcaption}

\definecolor{wheat}{rgb}{0.96,0.87,0.70}
\DeclareMathOperator*{\argmax}{argmax}
\DeclareMathOperator*{\argmin}{argmin}

\begin{document}
\begin{frontmatter}
\title{COBALT: COnstrained Bayesian optimizAtion of computationaLly expensive grey-box models exploiting derivaTive information}

\author[OSU]{Joel A. Paulson\corref{c1}}
\cortext[c1]{Corresponding author}
\ead{paulson.82@osu.edu}
\author[OSU]{Congwen Lu}
\ead{lu.2318@osu.edu}

\address[OSU]{Department of Chemical and Biomolecular Engineering, The Ohio State University, Columbus OH, USA}

\begin{abstract}
Many engineering problems involve the optimization of computationally expensive models for which derivative information is not readily available. The Bayesian optimization (BO) framework is a particularly promising approach for solving these problems, which uses Gaussian process (GP) models and an expected utility function to systematically tradeoff between exploitation and exploration of the design space. BO, however, is fundamentally limited by the black-box model assumption that does not take into account any underlying problem structure. In this paper, we propose a new algorithm, COBALT, for constrained \textit{grey-box} optimization problems that combines multivariate GP models with a novel constrained expected utility function whose structure can be exploited by state-of-the-art nonlinear programming solvers. COBALT is compared to traditional BO on seven test problems including the calibration of a genome-scale bioreactor model to experimental data. Overall, COBALT shows very promising performance on both unconstrained and constrained test problems.
\end{abstract}


\end{frontmatter}

%%%%%%%%%%%%%%%%%%%%%%%%%
\section{Introduction}
\label{sec:introduction}

Design problems, which can generally be formulated as mathematical optimization problems \cite{herskovits05}, occur in a wide-variety of science, engineering, and manufacturing endeavors. For example, pharmaceutical researchers must design new drugs to fight diseases, social media companies must design user-friendly websites to increase advertising revenue, and process engineers must synthesize flowsheets that achieve the desired goals of the process (e.g., profitable operation that meets chemical product specifications with minimal waste). In certain situations, one is able to develop an \textit{equation-oriented} (EO) model (also known as ``first-principles'', ``physics-based'' or ``white-box'' models) of the system whose structure can be exploited by existing solvers that take advantage of first- and/or second-order derivative information (see, e.g., \cite{sahinidis96, boggs00, biegler09, misener14}). However, obtaining accurate EO models for each and every component of a complex system is not always possible. Examples of non-EO models (also known as ``simulation-based'' or ``black-box'' models) include thermodynamic property relationships, models of a proprietary unit operations, and expensive finite-element, partial differential equation-based, and molecular simulations. 

When gradient information is not readily available, as is the case for black-box models, one of the main alternatives is to rely on so-called derivative-free optimization (DFO) methods \cite{conn09, rios13, larson19}, which can be very broadly divided into stochastic and deterministic approaches. The majority of the stochastic DFO methods can be classified as either evolutionary or population-based algorithms including genetic algorithms \cite{mukhopadhyay09}, particle swarm optimization (PSO) \cite{eberhart1995}, and the covariance matrix adaptation evolution strategy (CMA-ES) \cite{hansen2003reducing}. A key limitation of these methods, however, are that they often require a large number of function evaluations to find the optimum \cite{wessing17} such that they are not directly applicable to expensive black-box simulators.\footnote{Another important limitation of many DFO methods is that they cannot directly handle general nonlinear and/or black-box constraints and they instead focus on box constraints that can be easily handled using simple projection operators. Interested readers are referred to \cite[Section 7]{larson19} for more details on constraint handling methods in DFO.} Deterministic DFO methods, on the other hand, are often motivated by the optimization of an expensive objective function and can be classified as either \textit{direct search} or \textit{model-based} methods. Direct methods, which includes Nelder-Mead simplex algorithm \cite{nelder65}, mesh adaptive direct search (MADS) \cite{audet06}, and generalized pattern search (GPS) \cite{kolda03}, determine search directions from the evaluations of the objective function directly. Alternatively, the so-called model-based methods construct a surrogate model for the objective function using these evaluations to better guide the search process. Model-based DFO methods mainly differ by the choice of \textit{scale} (local versus global approximation) and the \textit{type} of function approximator (e.g., polynomial, neural network, or radial basis function models). 

Not only is selecting the ``right'' type of surrogate model challenging when little is known about the structure of the objective, this choice can have a strong effect on the performance of the model-based DFO algorithm. Gaussian process (GP) models are a particularly interesting class of surrogates due to the fact that they are \textit{probabilistic} and \textit{non-parametric}; GP models are easily derived by placing a prior over the set of possible objective functions that can be recursively updated to account for measured data (i.e., objective evaluations) using Bayes' rule \cite{rasmussen06}. By combining the GP model of the objective with an expected utility (or acquisition) function that leverages the uncertainty in the posterior distribution, we arrive at what is commonly referred to as the Bayesian optimization (BO) framework \cite{pelikan99, brochu10, shahriari15, frazier18}. A key ingredient in BO is the choice of the acquisition function that should be designed in a way that tradeoffs \textit{exploration} of regions where the surrogate model is uncertain and \textit{exploitation} of the model's confidence in good solutions. Although the basic BO framework can be traced back to the 1970s \cite{movckus75}, its popularity has substantially grown in recent years due to advances in computer power, algorithms, and software as well as successes in a variety of application areas including hyperparameter optimization in machine learning models \cite{bergstra11}, material design and discovery \cite{frazier16}, aircraft design \cite{meliani19}, and automated controller design \cite{sorourifar20,paulson20,sorourifar21}. 

Although BO was originally formulated for unconstrained problems, it has been recently extended to handle expensive-to-evaluate black-box constraints. There are two main classes of methods for handling these constraints, which we categorize as \textit{implicit} and \textit{explicit}. Implicit approaches define a new objective using a merit-type acquisition function that simultaneously accounts for effects due to the unknown objective and unknown constraints. Several merit functions have been proposed in the literature including the expected constrained improvement \cite{gardner14} and the augmented Lagrangian BO (ALBO) method \cite{picheny16} that combines the classical augmented Lagrangian method with unconstrained BO. Explicit approaches, on the other hand, attempt to model the constraints and solve a constrained subproblem that is restricted to a ``best guess'' of the feasible region, and can be further subdivided into \textit{deterministic} and \textit{probabilistic} methods. In deterministic explicit methods, one would disregard the variance information predicted by the GP model such as the super efficient global optimization (SEGO) approach in \cite{sasena02}. Probabilistic explicit methods instead embrace this uncertainty by either \textit{restricting} the feasible region to ensure high probability of constraint satisfaction at each iteration or, more recently, \textit{relaxing} the feasible region to allow exploration in the case of poorly modeled constraints. In particular, the recently proposed upper trust bound (UTB) method \cite{priem19}, which allows the mean prediction of the GP model to violate constraints up to a constant factor times the variance of the GP model, has shown the ability to more effectively compromise between exploration and exploitation of the feasible design space than alternative constrained BO methods.

Even though constrained BO methods have been found to empirically perform well on a variety of complex engineering problems in which the dimension of the design space is relatively small (typically on the order of ten or less), its sample efficiency tends to decrease with increasing dimension due to the exponential growth in the size of the search space. Although this challenge can partly be addressed by combining sensitivity and/or dimensionality reduction techniques with GP models (see, e.g., \cite{bouhlel16}), these black-box methods are fundamentally limited by the fact that any available knowledge about the structure of the underlying objective and constraint functions is neglected. In many practical engineering problems, only a portion (or subsystem) of the model is not explicitly known; such cases do not neatly fit into the either the white- or black-box problem classes and thus we can introduce the notion of a hybrid ``grey-box'' model that involves a mixture of EO and non-EO models. In this work, we represent grey-box models as \textit{composite} objective and constraint functions of the form $f(x) = g(h(x))$ where $g(\cdot)$ and $h(\cdot)$ are the white-box and black-box functions, respectively, which appear in many important real-world problems. For example, when calibrating parameters $x$ of an expensive process simulator to measured data $y_\text{meas}$, the objective function to be minimized can be formulated as $f(x) = g(h(x)) = \| h(x) - y_\text{meas} \|$ where $h(x)$ is the predicted output of the simulator for fixed parameters $x$ and $\| \cdot \|$ is some monotonic transformation of the likelihood of the measurement errors. It was recently shown in \cite{astudillo19} that significant improvements in the convergence rate of BO can be achieved when accounting for this composite structure for unconstrained problems. Extending this approach to constrained grey-box problems is not trivial, as it relies on a stochastic gradient ascent algorithm to maximize the composite acquisition function that is not directly applicable to nonlinear and non-convex constraints. 

A variety of methods for constrained grey-box optimization have been developed within the process systems engineering community \cite{eason16, eason18, bajaj18, boukouvala17, kim20, schweidtmann19, beykal18, beykal18b, beykal20}, which could be used as alternatives to the BO framework. One recent example is the trust region filter algorithm proposed in \cite{eason16, eason18}, which is guaranteed to converge to a local optimum. A potential disadvantage of this approach, however, is the lack of a global surrogate model, which may result in convergence to a highly suboptimal local optimum depending on the selected initialization point. The ARGONAUT algorithm \cite{boukouvala17} constructs a global surrogate model that is sequentially optimized using the global optimization solver ANTIGONE; a similar adaptive sampling framework has been recently extended to mixed-integer nonlinear programs in \cite{kim20}. There are two potential disadvantages to these (and many other related) methods: (i) they often do not directly account for the composite structure of the objective and/or constraint functions and (ii) only a deterministic surrogate model is trained, which tends to over exploit the initial runs whenever only a small number of function evaluations can be performed due to a limited computational budget. 
%JAP: left out alamo, which is bascially just use all of your buget to build model (lowest error adaptive sampling) and then optimize

Motivated by the BO framework and results in \cite{astudillo19}, we propose a novel algorithm for constrained grey-box optimization problems in this work, which we refer to as \textbf{COBALT} (\textbf{CO}nstrained \textbf{B}ayesian optimiz\textbf{A}tion of computationa\textbf{L}ly expensive grey-box models exploiting deriva\textbf{T}ive information). The proposed COBALT algorithm is composed of the following three main components: (i) a multivariate GP model of the black-box portions of the problem, (ii) a novel acquisition function for composite objective functions that is almost everywhere differentiable, and (iii) a generalization of the UTB constraint handling method to the case of composite functions using the notion of chance constraints. Due to the composite structure of the objective and constraints, we cannot derive the simple analytic expressions for the constrained acquisition function often found in traditional BO methods. Instead, we utilize the sample average approximation (SAA) method \cite{kleywegt02} to convert the stochastic constrained acquisition optimization problem into a deterministic problem. To alleviate the challenges that arise with an SAA-based reformulation of the chance constraints \cite{pagnoncelli09}, we propose a moment-based approximation that greatly simplifies the complexity of the SAA problem, which we show can be efficiently optimized using state-of-the-art NLP solvers. Through extensive testing on various types of test problems, we have observed that COBALT is able to outperform traditional BO by finding better (up to multiple orders of magnitude) quality solutions in fewer iterations, which can be directly attributed to its ability to account for the grey-box structure of the problem. 

The remainder of the paper is organized as follows. In Section \ref{sec:problem-definition}, the constrained grey-box problem of interest in this work is formulated. In Section \ref{sec:cobalt}, the proposed COBALT algorithm and its relevant parameters are presented. Section \ref{sec:numerical-example} discusses a Matlab-based implementation of COBALT and presents results and comparisons for six benchmark global optimization problems and a complex parameter estimation problem for a genome-scale bioreactor model. Lastly, we conclude the article and discuss some important directions for future work in Section \ref{sec:conclusions}. 

\subsection{Notation}

Throughout the paper, we use the following notation. We let $\| x \|_p$ denote the $\ell_p$ norm of a vector $x \in \mathbb{R}^n$. We let $\mathbb{S}_{+}^n$ and $\mathbb{S}_{++}^n$ denote the set of positive semidefinite and positive definite $n \times n$ matrices, respectively. 
%We denote vector concatenation by $(a,b) = [a^\top, b^\top]^\top \in \mathbb{R}^{n+m}$ for any vectors $a \in \mathbb{R}^n$ and $b \in \mathbb{R}^m$. 
By $\Pi_X$, we denote the Euclidean projection operator onto the set $X$, that is, $\Pi_X(x) = \argmin_{x' \in X} \| x - x' \|_2$. For any $a \in \mathbb{R}$, we denote $[a]^+ = \max\{ a, 0 \}$. The notation $\lfloor a \rfloor$ is the largest integer less than or equal to $a \in \mathbb{R}$ and $\lceil a \rceil$ is the smallest integer greater than or equal to $a \in \mathbb{R}$. For a real-valued function $f : \mathbb{R}^n \to \mathbb{R}$, we let $\nabla_x f(x) = (\partial f(x)/\partial x_1, \ldots, \partial f(x)/\partial x_n)$ denote its gradient with respect to $x$ and simplify this to $\nabla_x f(x) = \nabla f(x)$ when the argument is clear from the context. For random vector $X$, we let $\mathbb{E}_X \{ \cdot \}$ denote the expectation operator, $\mathbb{E}_{X | Y}\{ \cdot \}$ denote the conditional expectation given $Y$, and $\mathbb{P}\{ X \in A \}$ denote the probability that $X$ lies in the set $A$. By $\mathcal{N}(\mu, \Sigma)$, we denote a multivariate Gaussian distribution with mean $\mu \in \mathbb{R}^n$ and covariance $\Sigma \in \mathbb{S}_{+}^n$. The subscript $n$ exclusively refers to current iteration; additional subscripts may be added to denote elements and/or samples of particular variables.

%%%%%%%%%%%%%%%%%%%%%%%%%
\section{Problem Definition}
\label{sec:problem-definition}

In this paper, we consider a general constrained nonlinear grey-box optimization problem that can be mathematically defined as
\begin{subequations} \label{eq:grey-box-problem}
\begin{align}
\label{eq:grey-box-problem-f}
\min_{x, y, z} &~~ f(x, y), \\
\label{eq:grey-box-problem-g}
\text{s.t.} &~~ g(x, y) \leq 0, \\
\label{eq:grey-box-problem-d}
&~~ y = d(z), \\
\label{eq:grey-box-problem-A}
&~~ z = Ax, \\
\label{eq:grey-box-problem-xyz}
&~~ x \in \mathcal{X} \subset \mathbb{R}^{n_x}, ~ y \in \mathbb{R}^{n_y}, ~ z \in \mathbb{R}^{n_z}, 
\end{align}
\end{subequations}
where $x \in \mathbb{R}^{n_x}$ is a vector of $n_x$ decision variables that is constrained to a constrained within known lower $x^L$ and upper $x^U$ bounds, i.e., $\mathcal{X} = \{ x : x^L \leq x \leq x^U \}$; $z \in \mathbb{R}^{n_z}$ and $y \in \mathbb{R}^{n_y}$, respectively, denote the inputs and outputs of an unknown ``black-box'' vector-valued function $d : \mathbb{R}^{n_z} \to \mathbb{R}^{n_y}$ referenced in \eqref{eq:grey-box-problem-d}; $A \in \mathbb{R}^{n_x \times n_z}$ is a binary matrix that encodes that the black-box function may only require a subset $n_z \leq n_x$ of $x$ as input according to \eqref{eq:grey-box-problem-A}; and $f : \mathbb{R}^{n_x} \times \mathbb{R}^{n_y} \to \mathbb{R}$ and $g : \mathbb{R}^{n_x} \times \mathbb{R}^{n_y} \to \mathbb{R}^{n_g}$ are the known ``white-box'' objective and constraint functions, respectively. We emphasize the so-called ``grey-box'' structure of the optimization problem \eqref{eq:grey-box-problem} in that both $f$ and $g$ are \textit{composite} functions that are second-order continuously differentiable with known structure whereas $d$ is completely unknown. Here the phrase ``unknown'' strictly refers the mathematical expression of $d$ in terms of $z$; we do assume a simulator is available to query $y = d(z)$ at fixed $z \in \{ z : Ax^L \leq z \leq Ax^U \}$. Additionally, we assume this simulator is expensive to evaluate, so only a limited number of evaluations (on the order of tens to hundreds) can be performed.

Due to the vector representation of the expensive black-box function, we can easily incorporate any (finite) number of $n_b$ black-boxes by concatenating the outputs of each individual simulators
\begin{align}
d(Ax) = [d_1(A_1x)^\top , \ldots, d_{n_b}(A_{n_b} x)^\top ]^\top,
\end{align}
where $A_i \in \mathbb{R}^{n_x \times n_{z_i}}$ denotes the binary encoding matrix for the $i^\text{th}$ simulator for all $i = 1,\ldots, n_b$ and $A$ represents the collection of unique rows from $[A_1^\top, \ldots, A_{n_b}^\top]^\top$. 
%Furthermore, we can handle nonlinear equality constraints of the form $h(x,y) = 0$ by transforming them into two-sided inequalities $h(x, y) \leq 0$ and $h(x, y) \geq 0$ and incorporating them into the definition of $g$. 
%JAP: Change the equality constraints to a remark. 

An important distinction in the formulation of \eqref{eq:grey-box-problem} compared to alternative constrained grey-box algorithms is the explicit consideration of the composite structure of the overall objective $l(x) = f(x, d(Ax))$ and constraints $c(x) = g(x, d(Ax))$. Neglecting this structure, we obtain a simplified representation of \eqref{eq:grey-box-problem} as follows
\begin{subequations} \label{eq:simplified-grey-box}
\begin{align}
\min_{x \in \mathcal{X}} &~~ l(x), \\
\text{s.t.} &~~ c_u(x) \leq 0, ~~~ \forall u \in \{ 1,\ldots,U \}, \\
&~~ c_k(x) \leq 0, ~~~\forall k \in \{ 1, \ldots, K \},
\end{align}
\end{subequations}
where $l(x)$ must be generally modeled as an unknown function due to the embedded black-box function $d(Ax)$, $c(x) = [c_1(x), \ldots, c_{n_g}(x)]^\top$ is the vector concatenation of individual constraint functions, $k \in \{ 1, \ldots, K \}$ are the set of indices for which the constraints have known closed-form equations, $u \in \{ 1 ,\ldots, U\}$ are the set of indices for which the structure of the constraints are unknown, and $n_g  = K + U$. The formulation \eqref{eq:simplified-grey-box} has been considered in several previous works (see, e.g., \cite{boukouvala17, kim20, beykal18b, beykal20}) and is clearly a special case of \eqref{eq:grey-box-problem} whenever the composite structure of the unknown objective and constraints is neglected. Thus, one of the main contributions of this work is to exploit this composite structure (whenever possible) to improve the efficiency/performance of data-driven optimization procedures. The details of our proposed algorithm for solving \eqref{eq:grey-box-problem} are presented in the next section. 
% cite ARGONAUT, p-ARGOUNAUT, and Fani
% https://www.ncbi.nlm.nih.gov/pmc/articles/PMC6287910/#R8

%%%%%%%%%%%%%%%%%%%%%%%%%
\section{The COBALT Approach: Efficient Global Optimization of Constrained Grey-Box Models}
\label{sec:cobalt}

\begin{figure}[ht!]
\centering
\includegraphics[width=1.0\linewidth]{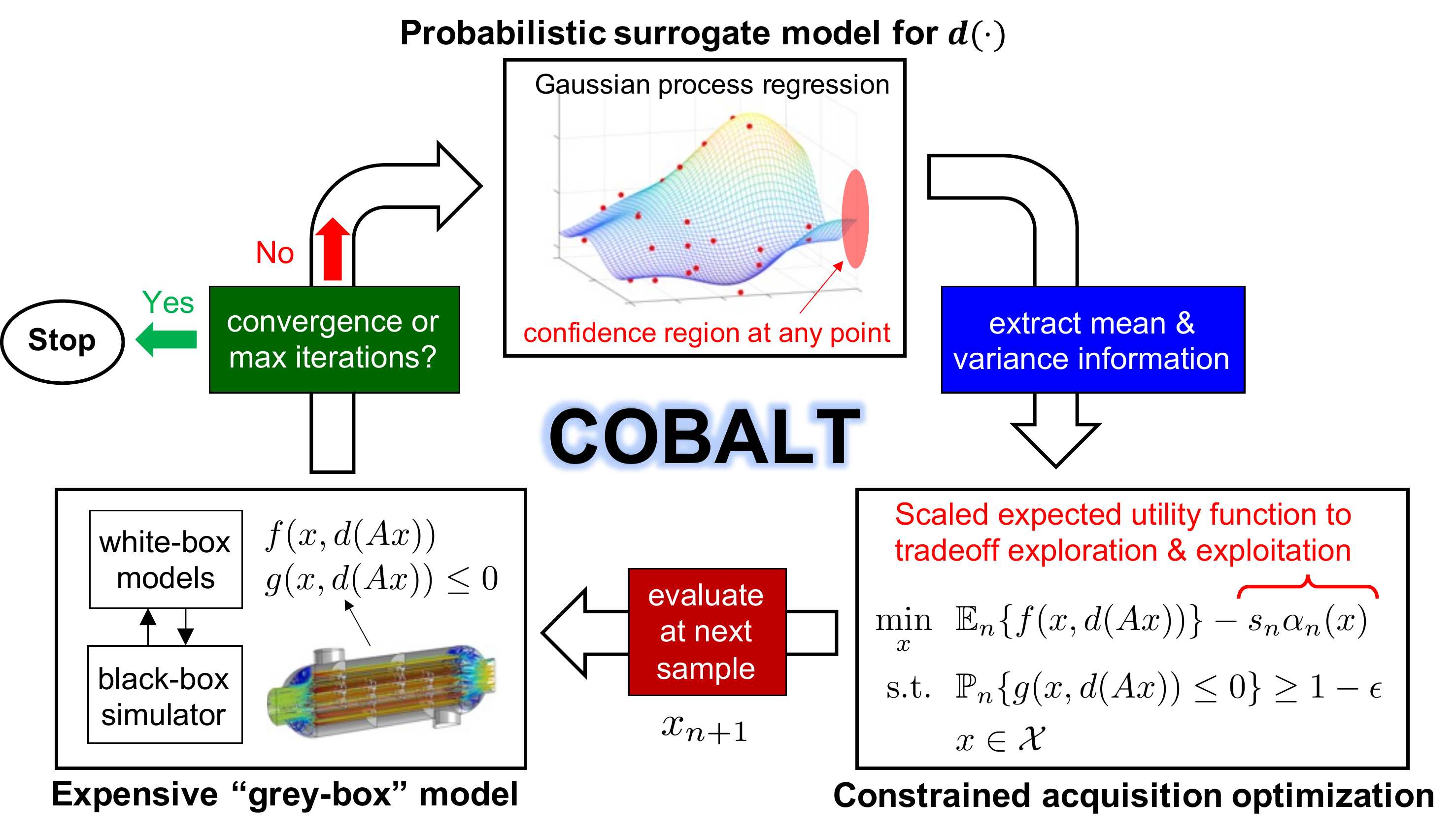}
\caption{Illustration of the main components of the proposed COBALT algorithm.}
\label{fig:cobalt_figure}
\end{figure}

The main idea behind the Bayesian optimization (BO) framework is to sequentially decide where to sample the design space $\mathcal{X}$ using all available observations of the objective and constraint functions. Instead of relying on measurements of $\ell(\cdot)$ and $c_u(\cdot)$, $\forall u \in \{ 1, \ldots, U \}$, as would be the case for traditional BO, we rely on observations of the black-box function $d$ in \eqref{eq:grey-box-problem} directly. Let $\mathcal{D}_n = \{ x_i, z_i, y_i \}_{i=1}^n$, where $z_i = Ax_i$ and $y_i = d(z_i)$, be all available observations at iteration $n$. By prescribing a prior belief over the function $d$, we can construct (and iteratively refine) a statistical surrogate model for $d$ given the available data using Bayes' rule to determine the posterior distribution $d(\cdot) | \mathcal{D}_n$. Given this probabilistic model, we can induce an \textit{acquisition} function $\alpha_n : \mathcal{X} \to \mathbb{R}$ that leverages uncertainty in the posterior to tradeoff between exploration and exploitation. We can roughly think of $\alpha_n$ as quantifying the \textit{utility} of potential candidate points for the next evaluation of $d$; by accounting for the composite structure we can focus sampling on regions that have the most significant impact on the overall objective $f(x, d(Ax))$ and likelihood of producing feasible points $g(x, d(Ax)) \leq 0$, as opposed to ``wasting'' evaluations in regions that are likely to yield poor results. A high-level overview of the BO process is provided in Algorithm \ref{alg:gp-sampling} and an illustrative flowchart is shown in Fig. \ref{fig:cobalt_figure}. 
%and illustrated in Fig. {\color{red} X}. 

\begin{algorithm}
\caption{General overview of the Bayesian optimization framework}
 \hspace*{\algorithmicindent} \textbf{Initialize:} Maximum number of samples $N$ and initial dataset $\mathcal{D}_1$. 
\begin{algorithmic}[1]
\For{$n=1$ to $N$}
\State Build statistical surrogate model for $d(\cdot)$;
\State Find $x_{n+1}$ as the solution to the following enrichment sub-problem
\begin{align} \label{eq:optimize-acqusition}
x_{n+1} = \argmax_{x \in \mathcal{X}_n} ~ \alpha_n(x),
\end{align}
\State Set $z_{n+1} = Ax_{n+1}$ and evaluate expensive function at $y_{n+1} = d(z_{n+1})$;
\State Augment the dataset $\mathcal{D}_{n+1}= \mathcal{D}_n \cup \{ x_{n+1}, z_{n+1}, y_{n+1} \}$.
\EndFor
\end{algorithmic}
 \hspace*{\algorithmicindent} \textbf{Output:} The feasible point with the lowest objective value. 
\label{alg:gp-sampling}
\end{algorithm}

There are three key ingredients in the proposed grey-box BO method summarized in Algorithm \ref{alg:gp-sampling}: (i) the choice of probabilistic surrogate model that consists of our prior beliefs about the behavior of $d$; (ii) the specification of acquisition functions $\alpha_n$ that captures the tradeoff between exploration and exploitation and can be relatively ``easily'' optimized; and (iii) the constraint handling mechanism by choice of the sets $\mathcal{X}_n \subseteq \mathcal{X}$ to ensure sufficient exploration of the feasible domain. The specific choices of these three elements that make up the proposed COBALT method are discussed in the remainder of this section.
%needed to define the enrichement subproblem \eqref{eq:optimize-acqusition}: 

\subsection{Overview of Gaussian process regression}
\label{subsec:GPsingle}
 
Gaussian process (GP) models represent an uncountable collection of random variables, any finite subset of which has a joint Gaussian distribution. Thus, GPs generalize the notion of a multivariate Gaussian distribution to ``distributions over functions,'' which are fully specified by their mean and covariance functions \cite{rasmussen06}. Although any probabilistic surrogate model can be used, e.g., \cite{snoek15}, we focus exclusively on GPs in this work due to their non-parametric nature, i.e., they can represent any function given a sufficiently large dataset. 

In this section, we provide an overview of GP regression for general scalar functions $s : \mathbb{R}^{n_z} \to \mathbb{R}$ from potentially noisy measurements
\begin{align} \label{eq:smeas}
t = s(z) + v,
\end{align}
where $v \sim \mathcal{N}(0, \sigma_v^2)$ is a zero mean Gaussian noise term with variance $\sigma_v^2$. Here, $s$ can be thought of as modeling a single component of the black-box function, i.e., $s = d_j$ for a given $j \in \{ 1,\ldots,n_y \}$; the extension of GP to multi-output functions will be discussed in more detail in the subsequent section. 

GPs are specified by their mean function $m(\cdot)$ and covariance function $k(\cdot, \cdot)$. We write that a function $s(\cdot)$ is distributed as a GP with mean function $m(\cdot)$ and covariance function $k(\cdot, \cdot)$ as follows
\begin{align} \label{eq:GPprior}
s(\cdot) \sim \mathcal{GP}(m(\cdot), k(\cdot, \cdot)),
\end{align}
with
\begin{subequations}
\begin{align}
m(z) &= \mathbb{E}_s\{ s(z) \}, \\
k(z,z') &= \mathbb{E}_s \{ (s(z) - m(z))(s(z') - m(z')) \},
\end{align}
\end{subequations}
where $z,z' \in \mathbb{R}^{n_z}$ are arbitrary input vectors and $\mathbb{E}_s\{ \cdot \}$ is the expectation over the function space. The GP prior can generally depend on a set of hyperparameters $\Psi_c$, i.e., $m(z | \Psi_c)$ and $k(z,z' | \Psi_c)$. Without loss of generality, we assume that the mean function is set to zero
\begin{align} \label{eq:zero-mean}
m(z | \Psi_c) = 0, 
\end{align}
which can be achieved by normalizing the data before training as discussed in, e.g., \cite{bradford18}. When using GP regression, the chosen class of covariance functions determines the properties of the fitted functions. Here, we will focus on stationary covariance functions from the Mat\'ern class whose smoothness can be adjusted by a parameter $\nu$ such that the corresponding function is $\lceil \nu/2 -1 \rceil$ times differentiable. Some of the most commonly used examples are
\begin{subequations} \label{eq:Matern-cov}
\begin{align}
k_{\nu=1}(z, z') &= \zeta^2 \exp\left(  -r(z,z') \right), \\
k_{\nu=3}(z, z') &= \zeta^2 \left(  1 + \sqrt{3}r(z,z') \right) \exp\left(  -\sqrt{3}r(z,z') \right), \\
k_{\nu=5}(z, z') &= \zeta^2 \left(  1 + \sqrt{5}r(z,z') + \frac{5}{3} r(z,z')^2 \right) \exp\left(  - \sqrt{5} r(z,z') \right),
\end{align}
\end{subequations}
where $r(z,z') = \sqrt{(z - z') \Lambda^{-2} (z - z') }$ is the scaled Euclidean distance and $\Lambda = \text{diag}(\lambda_1, \ldots, \lambda_{n_z})$ is a diagonal scaling matrix. Note that, in the limit $\nu \to \infty$, we recover the squared exponential (SQ-EXP) covariance function of the form
\begin{align} \label{eq:SE-cov}
k_{\text{SQ-EXP}}(z, z') &= \zeta^2 \exp\left(  -\frac{1}{2} r(z,z')^2 \right),
\end{align}
which is most commonly used function whenever $s$ is assumed to be a member of the space of smooth (infinitely differentiable) functions. 

Under assumption \eqref{eq:zero-mean}, the hyperparameters consist of $\Psi_c = [\zeta, \lambda_1, \ldots, \lambda_{n_z}]^\top$. The parameter $\zeta$ describes the output variance while the parameters $\{ \lambda_i \}_{i=1}^{n_z}$ define the length scale of each of the input variables. Covariance functions with different length scales for each input are called \textit{anisotropic}; if an input dimension is not important, then its corresponding length scale will be large. Due to the additive property of Gaussian distributions, we can derive a GP model for the observation
\begin{align}
t \sim \mathcal{GP}(0, k(z,z' | \Psi_c) + \sigma_v^2 \delta_{zz'}),
\end{align}
where $\delta_{zz'}$ is the Kronecker delta function that is equal to 1 whenever $z = z'$ and zero otherwise. If the function observations are noisy and $\sigma_v^2$ is unknown, it can be included in the joint set of hyperparameters for the prior denoted by $\Psi = [\Psi_c^\top, \sigma_v^2]^\top$. 

Training a GP model thus corresponds to calibrating the hyperparameters $\Psi$ to a given dataset. Let us assume that we have $n$ available measurements of the unknown function $s$ represented by the following matrices
\begin{subequations}
\begin{align}
Z &= [z_1, \ldots, z_n]^\top \in \mathbb{R}^{n \times n_z}, \\
T &= [t_1,\ldots, t_n]^\top \in \mathbb{R}^{n \times 1},
\end{align}
\end{subequations}
where $z_i$ and $t_i$ denote the $i^\text{th}$ input and output data point, respectively. Based on the GP prior assumption \eqref{eq:GPprior}, the measured data vector $T$ must follow a multivariate Gaussian distribution of the form
\begin{align}
T \sim \mathcal{N}(0, \Sigma_{T}), ~~~ [\Sigma_{T}]_{ij} = k(z_i, z_j | \Psi_c) + \sigma_v^2 \delta_{ij}, ~~~ \forall (i,j) \in \{ 1, \ldots, n \}^2.
\end{align}
We use the maximum likelihood estimation (MLE) framework to infer the hyperparameters $\Psi$ from the log-likelihood function of the observations
\begin{align}
\mathcal{L}(\Psi) = \log(p( T | Z, \Psi )) = -\frac{1}{2} T^\top \Sigma_{T}^{-1} T  - \frac{1}{2}\log(\det(\Sigma_{T})) - \frac{n}{2}\log(2\pi).
\end{align}
The MLE hyperparameter estimate is then specified as the solution to the following optimization problem
\begin{align} \label{eq:MLE}
\Psi_\text{MLE} = \argmax_{\Psi} ~  \mathcal{L}(\Psi),
\end{align}
which can be solved using readily available NLP methods. Once we have trained the hyperparameters, we can use the data $\mathcal{D} = \{ Z, T \}$ to infer the posterior distribution $s(z) | \mathcal{D}$ at any test point $z$ using Bayes' rule \cite[Chapter 2]{rasmussen06}
\begin{align} \label{eq:posterior-s}
s(z) | \mathcal{D} \sim \mathcal{N}(\mu_s(z ; \mathcal{D}), \sigma_s^2(z ; \mathcal{D})),
\end{align}
with
\begin{subequations} \label{eq:mean-variance-s}
\begin{align}
\mu_s(z ; \mathcal{D}) &= k(z, Z) \Sigma_{T}^{-1} T, \\
\sigma_s^2(z ; \mathcal{D}) &= k(z,z) - k(z,Z) \Sigma_{T}^{-1} k(z,Z)^\top, 
\end{align}
\end{subequations}
where
\begin{align}
k(z,Z) = [k(z,z_1), \ldots , k(z,z_n)] \in \mathbb{R}^{1 \times n}. 
\end{align} 
It is important to note that the posterior mean $\mu_s(z ; \mathcal{D})$ represents our best prediction of the unknown function $s(z)$ at any particular $z$ value, while the posterior variance $\sigma_s^2(z ; \mathcal{D})$ provides a measure of uncertainty in this prediction. 

\begin{remark} \label{rem:1}
The complexity of evaluating the posterior mean and variance scales as $O(n^3)$ with respect to number of observations $n$ due to inversion of the covariance matrix in \eqref{eq:mean-variance-s}. In practice, the Cholesky decomposition $\Sigma_T = L_T L_T^\top$ can be computed once and saved so that subsequent evaluations scale as $O(n^2)$, as long as the hyperparameters of the kernel are kept constant. This is typically not a major issue for expensive function evaluations (due to the relatively low computational budget which is the case of interest in this paper); however, does become an important challenge large datasets. There have been a significant number of contributions on reducing the computational cost including sparse GP methods \cite{snelson06}. This topic remains a very active area of research and recently developed packages such as GPyTorch \cite{gardner18} have been able to scale GPs to $n > 10^6$ training points. 
\end{remark}
% Liu, H., Cai, J., and Ong, Y.-S. Remarks on multi-output gaussian process regression. Knowledge-Based Systems, 144:102–121, 2018.
% Pleiss, G., Gardner, J. R., Weinberger, K. Q., and Wilson, A. G. Constant-time predictive distributions for gaussian processes. arXiv preprint arXiv:1803.06058, 2018.

\subsection{Statistical model for multi-output black-box function}

In traditional BO methods, a separate GP would be trained for the objective and unknown constraints in \eqref{eq:simplified-grey-box} following the procedure discussed in the previous section. Here, our goal is to learn a GP model for $d$ in \eqref{eq:grey-box-problem} instead using the dataset $\mathcal{D}_n$ that is recursively updated according to Algorithm \ref{alg:gp-sampling}. We model $d$ as being drawn from a multi-output GP distribution, i.e., $d(\cdot) \sim \mathcal{GP}(m(\cdot), k(\cdot, \cdot))$ where $m : \mathbb{R}^{n_z} \to \mathbb{R}^{n_y}$ is the prior mean function and $k : \mathbb{R}^{n_z} \times \mathbb{R}^{n_z} \to \mathbb{S}_{++}^{n_y}$ is the prior covariance function. Similarly to the single output case described previously, the posterior distribution $d(\cdot) | \mathcal{D}_n$ is again a multi-output GP (MOGP), $\mathcal{GP}(\mu_n(\cdot), K_n(\cdot, \cdot))$, where the posterior mean $\mu_n$ and covariance $K_n$ can be computed in closed-form \cite{liu18}. 

When modeling the correlation between the components of $d$, the evaluation cost of the posterior multi-output GP would scale as $O(n_y^2 n^3)$ in the worst-case (see Remark \ref{rem:1}). An alternative approach that we pursue here is to model the components of $d$ independently, meaning we constrain $K_n$ to be a diagonal matrix, so the necessary computations scale linearly with respect to the number of outputs $O(n_y n^3)$ in the worst-case. We focus on the more tractable case in this work since we need to systematically optimize over the MOGP model embedded within the acquisition and constraint functions, and we would like to limit the complexity of these sub-problems. Note that the proposed COBALT method, discussed in detail below, can flexibly handle \textit{any} MOGP model of interest and will provide immediate gains in performance when a more accurate MOGP model (i.e., one that better captures the underlying correlation between the elements of $d$) is utilized. However, as shown in \cite{liu18}, correlated MOGPs have a larger number of hyperparameters that must be estimated during the training procedure, which makes the MLE estimation problem more difficult to solve. Therefore, in practice, correlated MOGPs are not guaranteed have a higher prediction quality than uncorrelated MOGPs, especially in the low-data regime of interest in this work.

%Let $d(z) = [d_1(z), \ldots, d_{n_y}(z)]^\top$ %and $y_i = d_i(z)$ be the $i^\text{th}$ component of the black-box function for all $i \in \{ 1, \ldots n_y \}$. 
To build a separate GP for each $d_i(z)$ for all $i \in \{ 1, \ldots n_y \}$, we divide the complete dataset $\mathcal{D}_n = \{ \mathcal{D}_{1,n}, \ldots, \mathcal{D}_{n_y,n} \}$ into its individual components $\mathcal{D}_{i,n} = \{ x_{j}, z_j, d_i(z_j) \}_{j=1}^n$. Using the procedure summarized for scalar functions in Section \ref{subsec:GPsingle}, the posterior Gaussian distribution of $d(\cdot)$ at any test input $z$ is then
\begin{align} \label{eq:dGP}
d(z) | \mathcal{D}_n \sim \mathcal{N}(\mu_n(z), \Sigma_n(z)), 
\end{align}
with 
\begin{subequations} \label{eq:dmeanvariance}
\begin{align}
\mu_n(z) &= [\mu_{d_1}(z ; \mathcal{D}_{1,n}), \ldots, \mu_{d_{n_y}}(z ; \mathcal{D}_{n_y, n})]^\top, \\
\Sigma_n(z) &=  K_n(z, z) = \text{diag}( \sigma_{d_1}^2(z ; \mathcal{D}_{1,n}), \ldots, \sigma_{d_{n_y}}^2(z ; \mathcal{D}_{n_y,n}) ),
\end{align}
\end{subequations}
where $\mu_{d_i}(z ; \mathcal{D}_{i,n})$ and $\sigma^2_{d_i}(z ; \mathcal{D}_{i,n})$ are the posterior mean and variance functions for $d_i(z)$, respectively, built from the datasets $\mathcal{D}_{i,n}$ for all $i = 1,\ldots, n_y$.

\subsection{Modified expected improvement for composite functions}

Now that we have a statistical model to represent our belief about the unknown function $d$ at iteration $n$ as shown in \eqref{eq:dGP}, we need to select an acquisition function that captures the utility of sampling at a subsequent point $x_{n+1}$. We first focus on the unconstrained case for simplicity; the developed approach is extended to handle constraints in the next section. 

In experimental design and decision theory literature, the function $\alpha_n$ is often referred to as the \textit{expected utility}; whereas, in the BO literature, it is often called the \textit{acquisition} or \textit{infill} function. The acquisition function must be chosen carefully to achieve a reasonable tradeoff between exploring the search space and exploiting currently known promising areas of $\mathcal{X}$. As discussed in \cite{shahriari15}, acquisition functions can be categorized as either improvement-based, information-based, or optimistic. We develop a modified improvement-based policy in this work based on an extension of classical expected improvement (EI) to composite functions \cite{astudillo19}:
\begin{align} \label{eq:ei-cf}
\text{EI-CF}_n(x) = \mathbb{E}_n \left\lbrace [\ell_n^\star - f(x, d(Ax))]^+ \right\rbrace, 
\end{align}
where $\ell_n^\star = \min_{i \in \{ 1,\ldots, n \}} \ell(x_i)$ is the minimum value across the points that have been evaluated so far (often referred to as the \textit{incumbent}) and $\mathbb{E}_n\{ \cdot \}$ is the expected value conditioned on all available observations $\mathcal{D}_n$. 
%and $[a]^+ = \max(0, a)$ is the positive part function. 
When $d$ is scalar-valued ($n_y=1$), $A = I_{n_x}$, and $f$ is the identify function (i.e., $f(x,d(Ax)) = d(x)$), then \eqref{eq:ei-cf} reduces to the traditional EI function that can be computed analytically as follows
\begin{align}
\text{EI}_n(x) = (\ell_n^\star - \mu_n(x)) \Phi \left(  \frac{\ell_n^\star - \mu_n(x)}{\sigma_n(x)} \right) + \sigma_n(x)\phi \left(  \frac{\ell_n^\star - \mu_n(x)}{\sigma_n(x)} \right),
\end{align}
where $\Phi$ and $\phi$ are the standard Gaussian cumulative and probability density functions, respectively. Although such an analytic expression is relatively easy to evaluate and optimize, we are unable to derive one for $\text{EI-CF}_n(x)$ when $f$ is nonlinear in its second argument. Despite this potential complication, it can be shown that $\text{EI-CF}_n(x)$ is differentiable almost everywhere under mild regularity conditions, which are summarized next. We first recognize that, for any fixed $x \in \mathcal{X}$, the posterior distribution of $d(Ax)$ is a multivariate normal according to \eqref{eq:dmeanvariance}. We can thus exploit the ``whitening transformation'' to derive
\begin{align} \label{eq:whiteningtransform}
d(Ax) | \mathcal{D}_n = \mu_n(Ax) + C_n(Ax) \xi, 
\end{align}
where $\xi \sim \mathcal{N}(0,I_{n_y})$ is a standard normal random vector and $C_n(z)$ is the lower Cholesky factor of $\Sigma_n(z)$. Substituting this transformation into \eqref{eq:ei-cf}, we can replace $\mathbb{E}_n\{ \cdot \}$ with an expectation over $\xi$, i.e., $\mathbb{E}_\xi\{ \cdot \}$, implying that $\text{EI-CF}_n(x)$ can be straightforwardly estimated using Monte Carlo (MC) sampling; note $\text{EI-CF}_n(x)$ is finite for all $x \in \mathcal{X}$ whenever $\mathbb{E}_\xi\{ | f(x, \xi) | \} < \infty$. Similarly to \cite[Supplementary Material]{astudillo19}, we now make the following mild assumptions about $f$, $\mu_n$, and $\Sigma_n$:
\begin{assumption} \label{assumption:1}
Let $\mathcal{X}_0$ be an open subset of $\mathcal{X}$ so that $\mu_n(Ax)$ and $\Sigma_n(Ax)$ are differentiable for every $x \in \mathcal{X}_0$. A measurable function $\eta : \mathbb{R}^{n_y} \to \mathbb{R}$ exists such that
\begin{enumerate}[1.]
\item The function $f$ is differentiable;
\item $\| \nabla f(x, \mu_n(Ax) + C_n(Ax) \xi) \|_2 < \eta(\xi)$ for all $x \in \mathcal{X}_0$ and $\xi \in \mathbb{R}^{n_y}$;
\item $\mathbb{E}_\xi\{ \eta(\xi)^2 \} < \infty$ is finite for $\xi \sim \mathcal{N}(0, I_{n_y})$;
\item $\{ x \in \mathcal{X}_0 : f(x, \mu_n(Ax) + C_n(Ax) \xi) = \ell_n^\star \}$ is countable for almost every $\xi \in \mathbb{R}^{n_y}$.
\end{enumerate}
\end{assumption}

As long as the prior mean function $m(Ax)$ and covariance function $K(Ax,Ax)$ are differentiable on $\text{int}(\mathcal{X})$, we can show that $\mu_n(Ax)$ and $\Sigma_n(Ax)$ are differentiable on $\mathcal{X}_0 = \text{int}(\mathcal{X}) \setminus \{ x_1, \ldots, x_n \}$. Combining this with Assumption \ref{assumption:1}, we see $\text{EI-CF}_n(x)$ must be differentiable almost everywhere (except on a countable subset of $\mathcal{X}$). Thus, when it exists, the gradient of the composite EI function is given by
\begin{align} \label{eq:ei-cf-gradient}
\nabla \text{EI-CF}_n(x) = \mathbb{E}_\xi\{ \gamma_n(x, \xi) \},
\end{align}
where 
\begin{align} \label{eq:gamma-n}
\gamma_n(x,\xi) = \begin{cases}
-\nabla f(x, \mu_n(Ax) + C_n(Ax)\xi), &\text{if} ~ f(\mu_n(Ax) + C_n(Ax) \xi) < f_n^\star, \\
0, &\text{otherwise}.
\end{cases}
\end{align}
In the derivation of \eqref{eq:ei-cf-gradient}, we have switched the order of the expectation and the gradient operators that is generally possible under minor technical conditions presented in \cite[Theorem 1]{lecuyer90}. Since $\gamma_n(x, \xi)$ is an unbiased estimator of $\nabla \text{EI-CF}_n(x)$, we could potentially use the following stochastic gradient ascent algorithm to search for maximizer of \eqref{eq:optimize-acqusition}, with $\alpha_n \leftarrow \text{EI-CF}_n$, from some initial guess $x_{n+1, 0}$:
\begin{align} \label{eq:sga}
x_{n+1, t+1} = \Pi_{\mathcal{X}_{n}} \left( x_{n+1, t} + \nu_t \gamma_n(x_{n+1, t}, \xi_t) \right), ~~ \forall t \in \{ 0, \ldots, T-1 \},
\end{align}
where $\{ x_{n+1,t} \}_{t=0}^{T}$ is the sequence of design variables over $T$ iterations of the algorithm,
%$\Pi_\mathcal{C}(\cdot)$ is the projection operator that projects its argument onto the set $\mathcal{C}$, 
$\{ \nu_t\}_{t=0}^{T-1}$ is the sequence of step sizes, and $\{ \xi_t \}_{t=0}^{T-1}$ are independent and identically distributed (i.i.d.) samples from the distribution of $\xi$. Even in the box-constrained case, i.e., $\mathcal{X}_{n} = \mathcal{X}$ such that the projection operator is simple to implement, the stochastic gradient ascent algorithm is known to be quite sensitive to the choice of step sizes $\nu_0, \ldots, \nu_{T-1}$ \cite{huan14}. Another key challenge is that $\gamma_n(x, \xi_i) = 0$ for a potentially large range of $x$ and $\xi$ values, as seen in \eqref{eq:gamma-n}, which may cause the iteration scheme \eqref{eq:sga} to become stuck locally. This particular problem only gets worse as the number of observations $n$ increase due to the fact that the incumbent $\ell_n^\star$ can only improve (or stay the same) at each iteration resulting in smaller probability that the surrogate predicts potential improvement. 
% Gradient-based stochastic optimization methods in Bayesian experimental design Xun Huan, Youssef M. Marzouk

This challenge has been previously noticed in the context of EI, which is often observed to be highly multimodal in nature \cite{sasena02}. One potential remedy is to modify the definition of the acquisition function to better reflect the shape of the objective. One such example is the ``locating the regional extreme'' acquisition proposed by Watson and Barnes in \cite{watson95} for black-box optimization problems, i.e., $\text{WB2}_n(x) = - \hat{\ell}_n(x) + \text{EI}_n(x)$ where $\hat{\ell}_n(x)$ is the mean of a posterior GP surrogate model for the overall objective. Although producing a smoother function, the $\text{WB2}_n(x)$ acquisition function does not account for the difference in the scales of the predicted objective and EI. In particular, as the GP surrogate model becomes more accurate, we expect $\text{EI}_n(x)$ to steadily decrease, resulting in $\text{WB2}_n(x) \approx - \hat{\ell}_n(x)$, so that Algorithm \ref{alg:gp-sampling} begins to fully focus on \textit{exploitation} of the surrogate, which may be undesired. Therefore, not only are we interested in improving the scaling of the WB2 function, we want to extend it to composite functions as follows
\begin{align} \label{eq:modified-WB2CF}
& \text{mWB2-CF}_n(x) = s_n \text{EI-CF}_n(x) - \hat{\ell}_n(x),
\end{align}
%&= \mathbb{E}_\xi \left\lbrace s_n [\ell_n^\star - f(x, \mu_n(Ax) + C_n(Ax) \xi)]^+ - f(x, \mu_n(Ax) + C_n(Ax) \xi) \right\rbrace,
where $s_n > 0$ denotes a non-negative scaling factor and $\hat{\ell}_n(x) = \mathbb{E}_n\{ f(x, d(Ax)) \}$ denotes the predicted mean of the overall objective function. 
%we have substituted the transformation in \eqref{eq:whiteningtransform} to derive the second line. 
The gradient can be straightforwardly computed using \eqref{eq:whiteningtransform}, \eqref{eq:ei-cf-gradient}, and \eqref{eq:gamma-n} as follows
\begin{align}
\nabla \text{mWB2-CF}_n(x) = \mathbb{E}_\xi\{ s_n \gamma_n(x, \xi) - \nabla f(x, \mu_n(Ax) + C_n(Ax) \xi) \},
\end{align}
which still exists almost everywhere and whose estimates does not suffer from the same zero gradient issue as $\gamma_n(x, \xi)$. This means we could use the same stochastic gradient ascent algorithm to solve the enrichment sub-problem \eqref{eq:optimize-acqusition} with the acquisition function set to $\alpha_n \leftarrow \text{mWB2-CF}_n$ by replacing the gradient estimate with $\gamma_n(x, \xi_i) \leftarrow s_n \gamma_n(x, \xi_i) - \nabla f(x, \mu_n(Ax) + C_n(Ax) \xi_i)$ in \eqref{eq:sga}. The only remaining question is how to select the scaling factor $s_n$ in such a way that approximately preserves the global optimum of the composite EI function
\begin{align}
\argmax_{x \in \mathcal{X}_{n}} ~ \text{EI-CF}_n(x) \approx \argmax_{x \in \mathcal{X}_{n}}  ~\text{mWB2-CF}_n(x).
\end{align}
Enforcing this condition is difficult in practice since we must account for the full variation of $\hat{l}_n(x)$ over $x \in \mathcal{X}_n$. Motivated by recent work in \cite{bartoli19}, we instead rely on a heuristic approach that uses a finite number of starting points $\mathcal{X}_{n, \text{start}} \subset \mathcal{X}_n$ to approximately maximize the composite EI, i.e., $\hat{x}_{n+1} \approx \argmax_{x \in \mathcal{X}_{n, \text{start}}} \text{EI-CF}_n(x)$ where $\mathcal{X}_{n, \text{start}}$ contains the starting points to be used in a multistart version of \eqref{eq:sga} (i.e.,  multiple restarts from different quasi-randomly sampled initial conditions $x_{n+1, 0}$). Using this approximate maximizer, we compute the scaling factor as
\begin{align}
s_n = \begin{cases}
\beta \frac{| \hat{\ell}_n(\hat{x}_{n+1}) |}{\text{EI-CF}_n(\hat{x}_{n+1})}, &\text{if} ~ \text{EI-CF}_n(\hat{x}_{n+1}) > 0, \\
1, &\text{otherwise},
\end{cases}
\end{align}
where $\beta > 1$ is an additional scaling factor that accounts for the degree of nonlinearity in $\hat{\ell}_n(x)$. We found that a relatively large value of $\beta = 100$ gives good results for a variety of functions studied in Section \ref{sec:numerical-example}. Finally, we highlight that, since neither $\text{EI-CF}_n(x)$ or $\hat{\ell}_n(x)$ can be exactly computed due to presence of nonlinear terms, these terms must be estimated by MC sampling. It is important to recognize, however, that the transformation \eqref{eq:whiteningtransform} implies that the relatively expensive evaluation of $\mu_n(Ax)$ and $C_n(Ax)$ only need to be done once at a given $x \in \mathcal{X}$ such that these MC estimates can be efficiently computed for even a relatively large number of samples of $\xi$. Furthermore, we can easily replace MC sampling with more computationally efficient uncertainty propagation methods, such as unscented transform \cite{julier97} or polynomial chaos expansions \cite{paulson17}, if and when needed. 

\subsubsection{Illustrative example comparing black- and grey-box formulations}

\begin{figure}[ht!]
        \centering
        \begin{subfigure}[b]{0.31\textwidth}
            \centering
            \includegraphics[scale=.5, width=\textwidth]{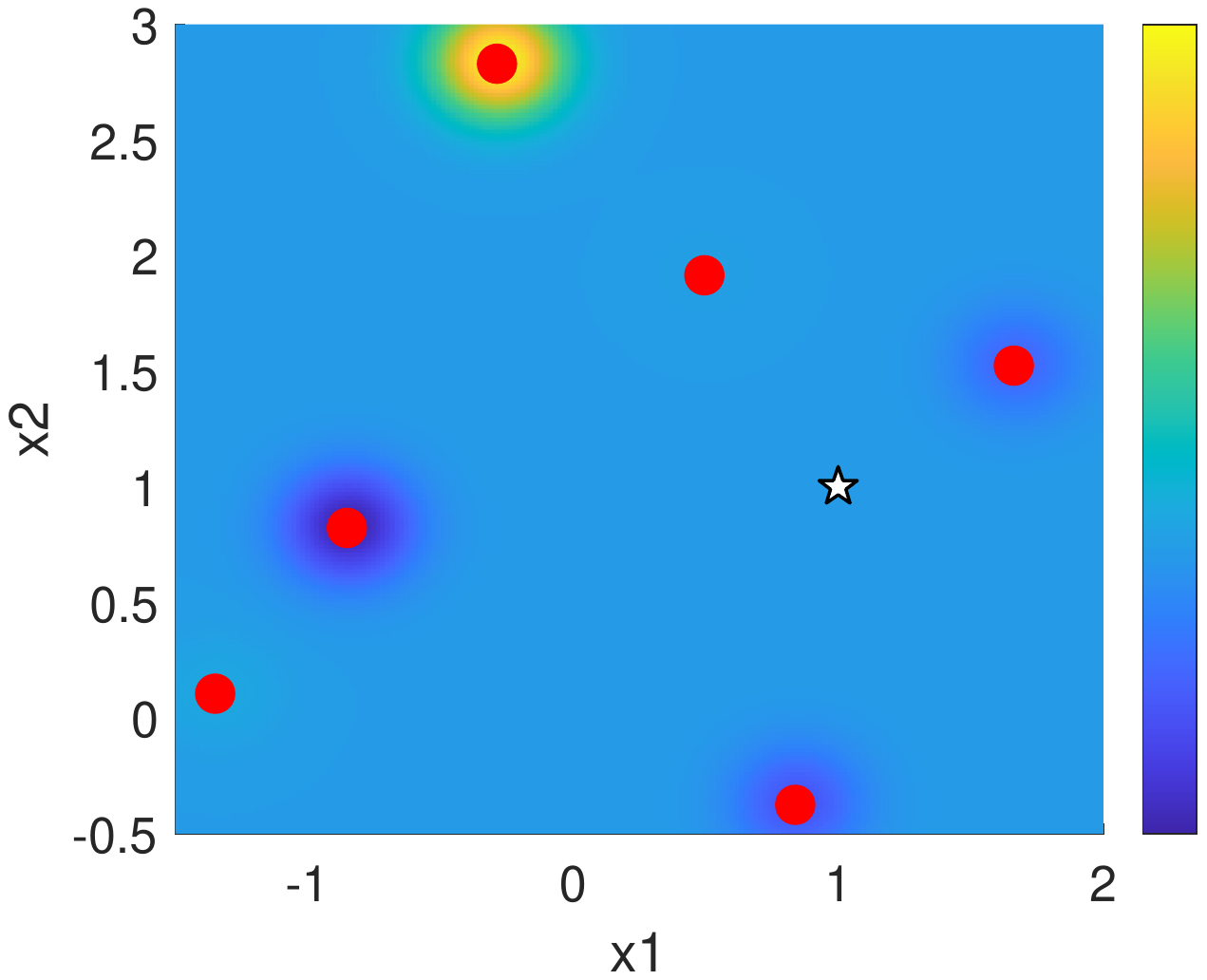}
            \caption{}
        \end{subfigure}
        \hfill
        \begin{subfigure}[b]{0.31\textwidth}  
            \centering 
            \includegraphics[scale=.5, width=\textwidth]{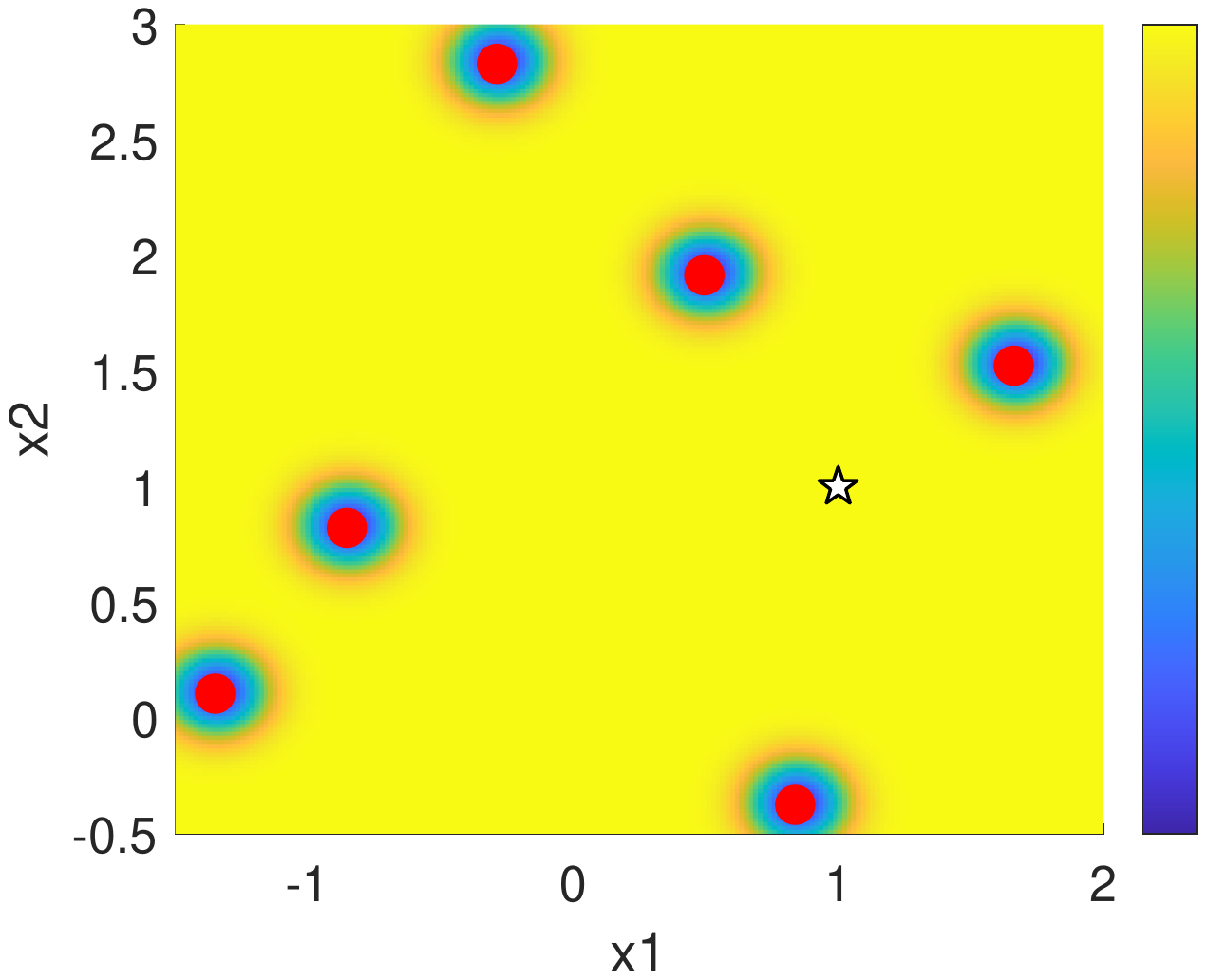}
            \caption{}
        \end{subfigure}
        \hfill
        \begin{subfigure}[b]{0.31\textwidth}   
            \centering 
            \includegraphics[scale=.5, width=\textwidth]{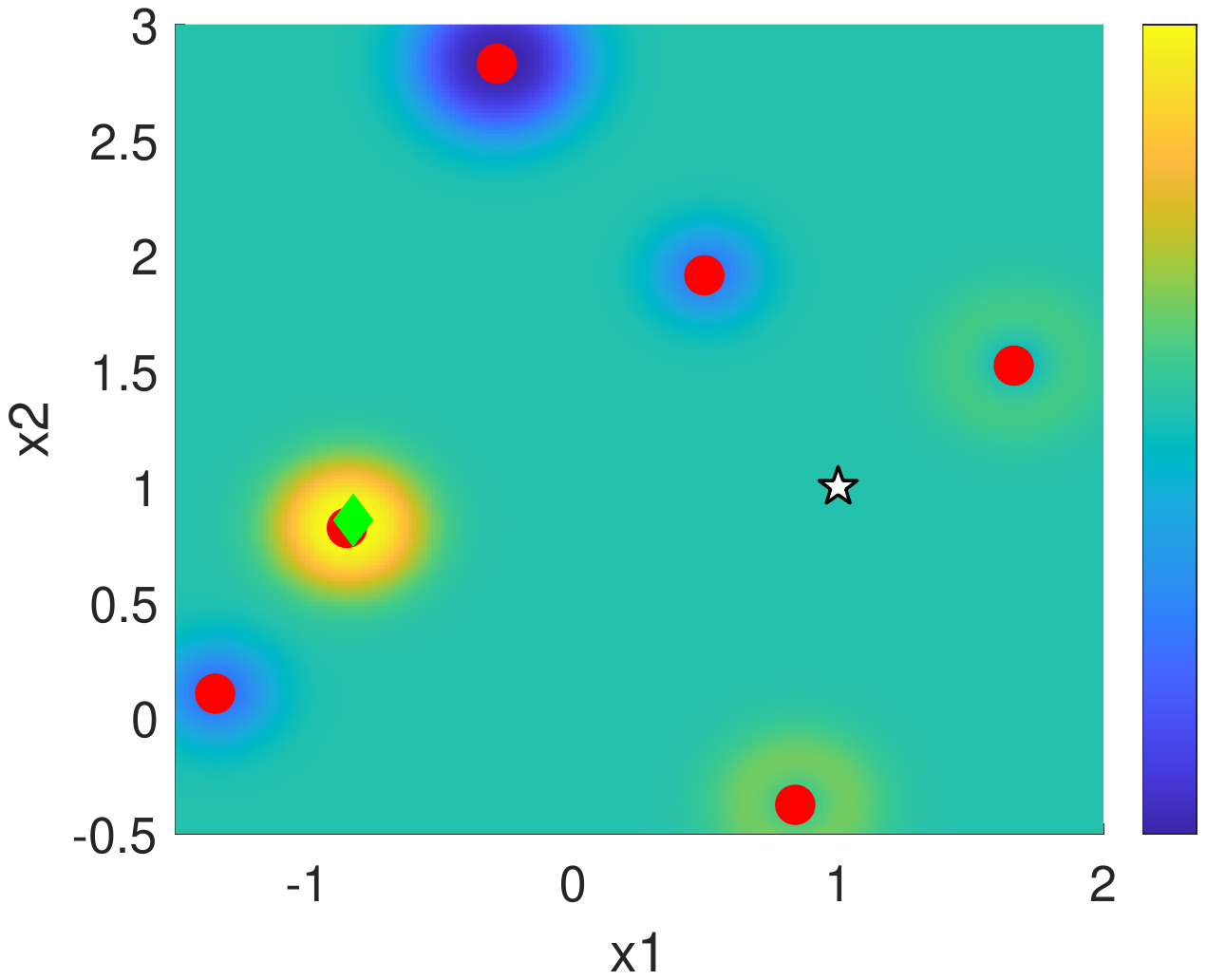}
            \caption{}
        \end{subfigure}
        \vskip\baselineskip
        \begin{subfigure}[b]{0.31\textwidth}
            \centering
            \includegraphics[scale=.5, width=\textwidth]{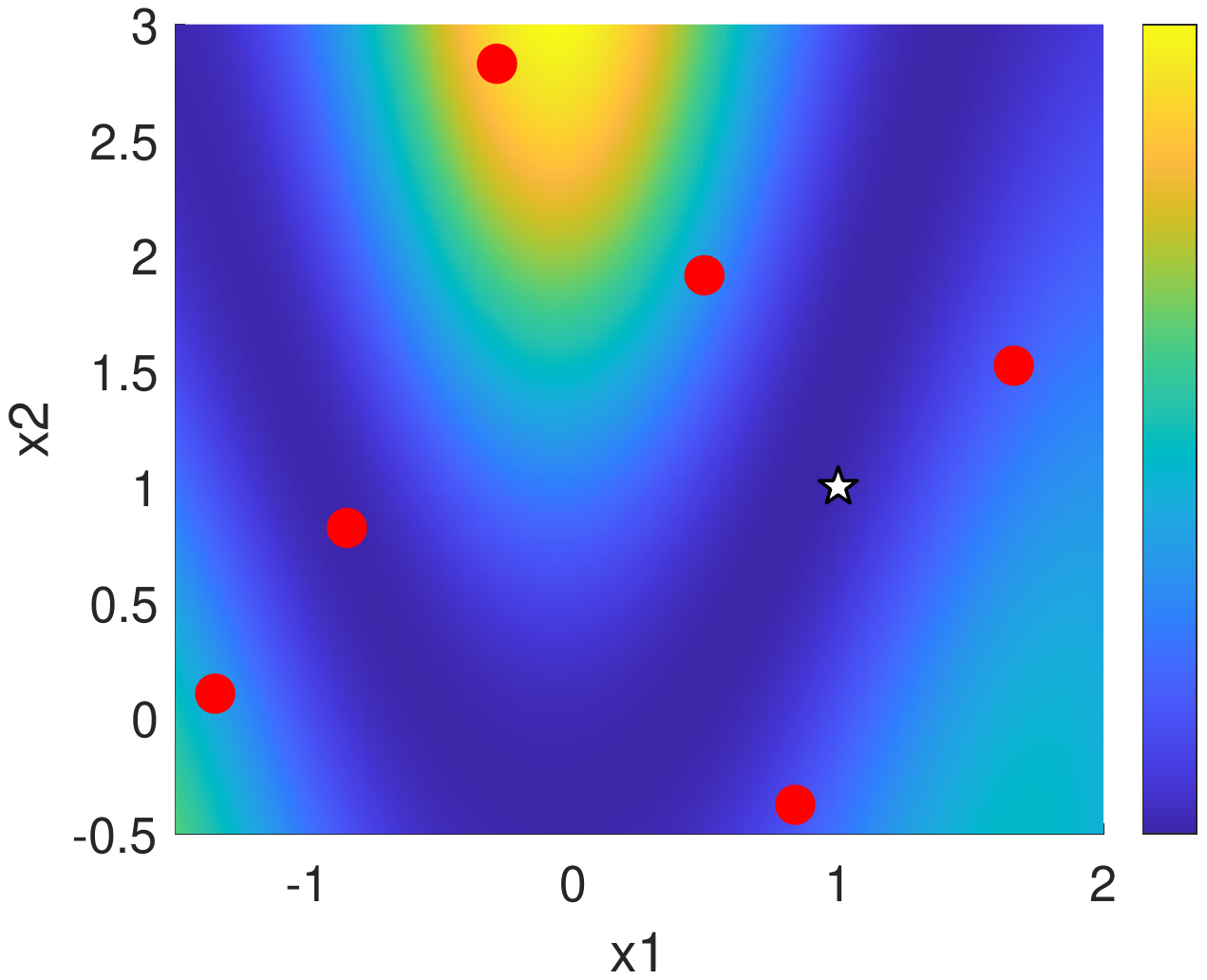}
            \caption{}
        \end{subfigure}
        \hfill
        \begin{subfigure}[b]{0.31\textwidth}
            \centering
            \includegraphics[scale=.5, width=\textwidth]{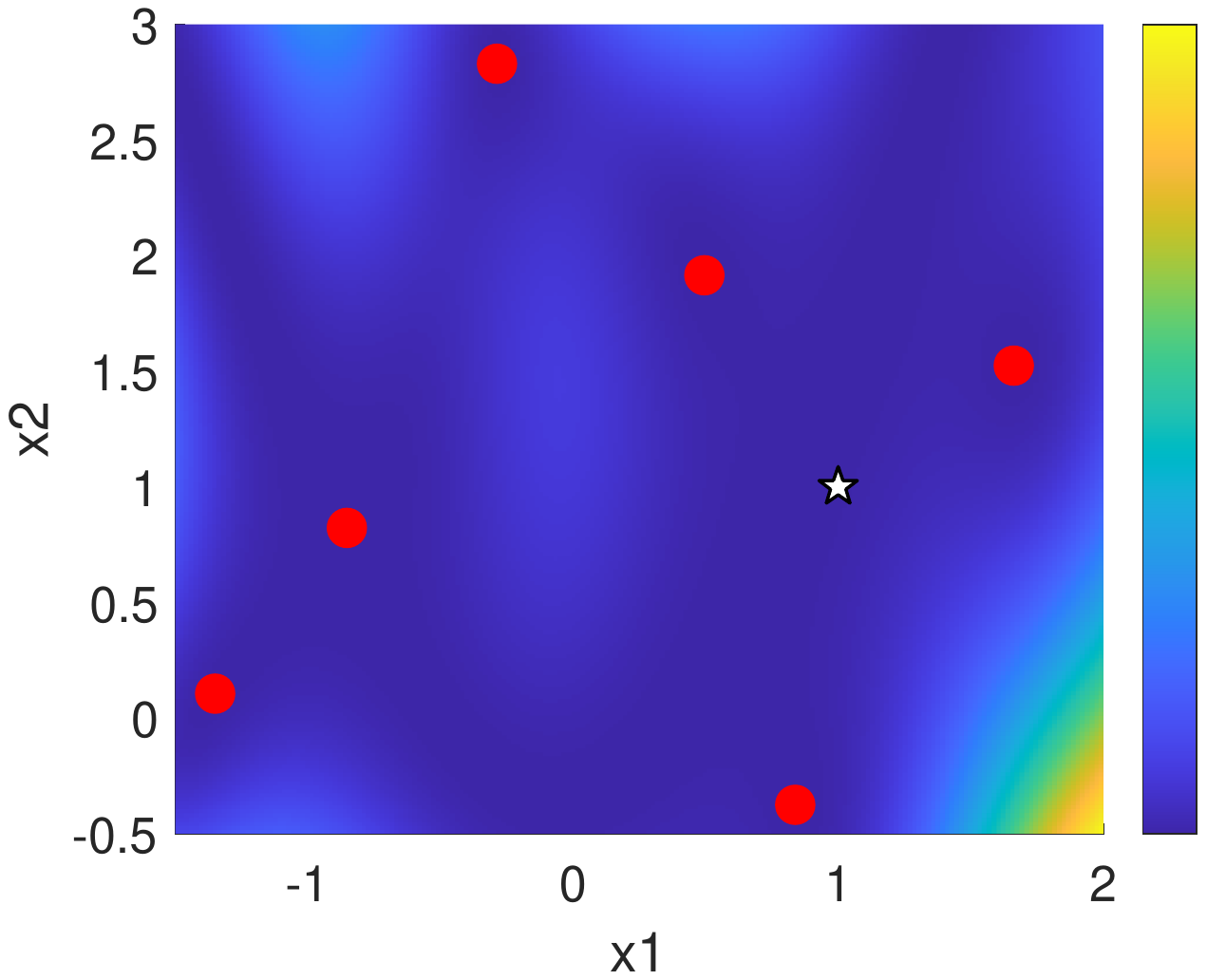}
            \caption{}
        \end{subfigure}
        \hfill
        \begin{subfigure}[b]{0.31\textwidth}
            \centering
            \includegraphics[scale=.5, width=\textwidth]{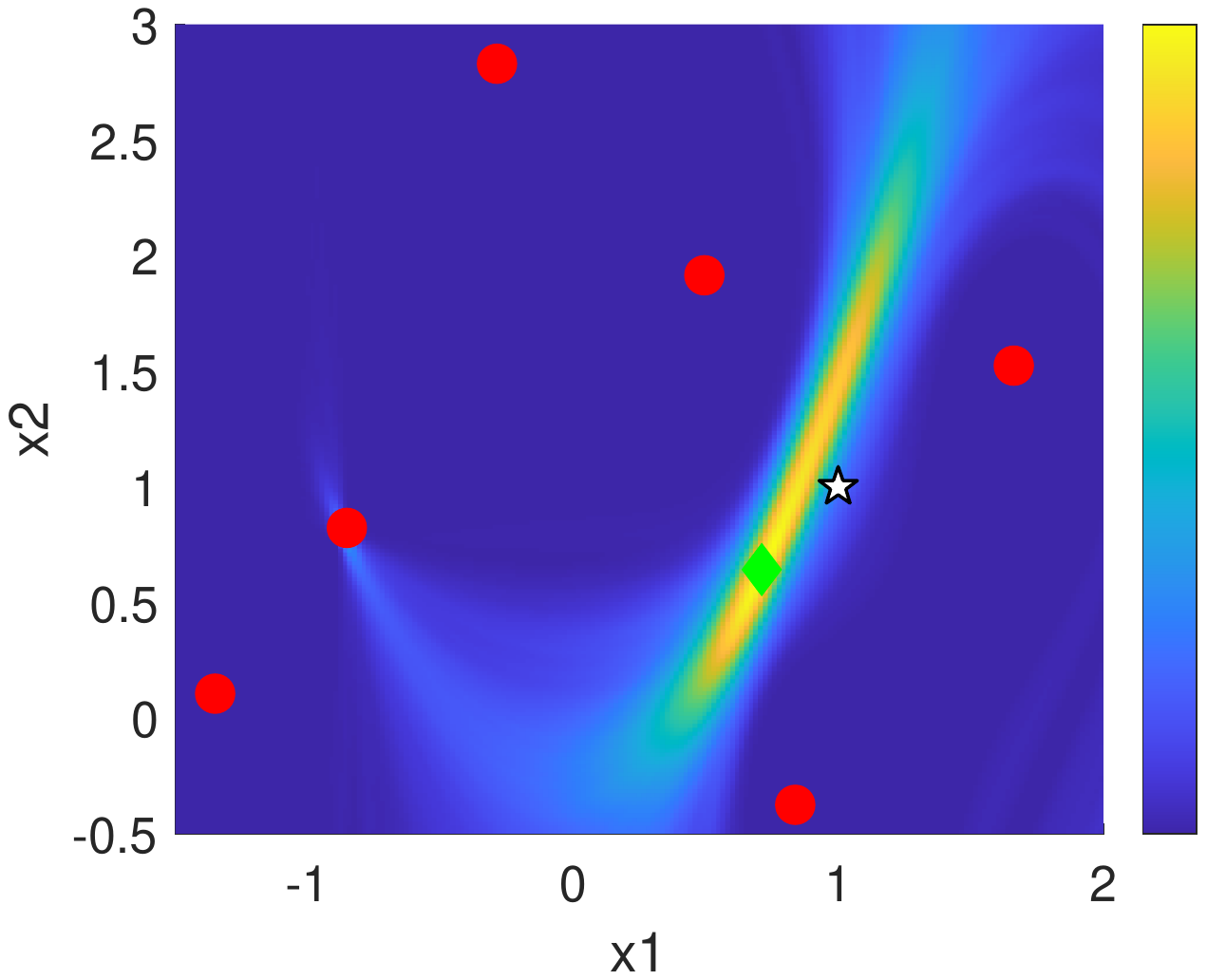}
            \caption{}
        \end{subfigure}
        \vskip\baselineskip
        \begin{subfigure}[b]{0.31\textwidth}
            \centering
            \includegraphics[scale=.5, width=\textwidth]{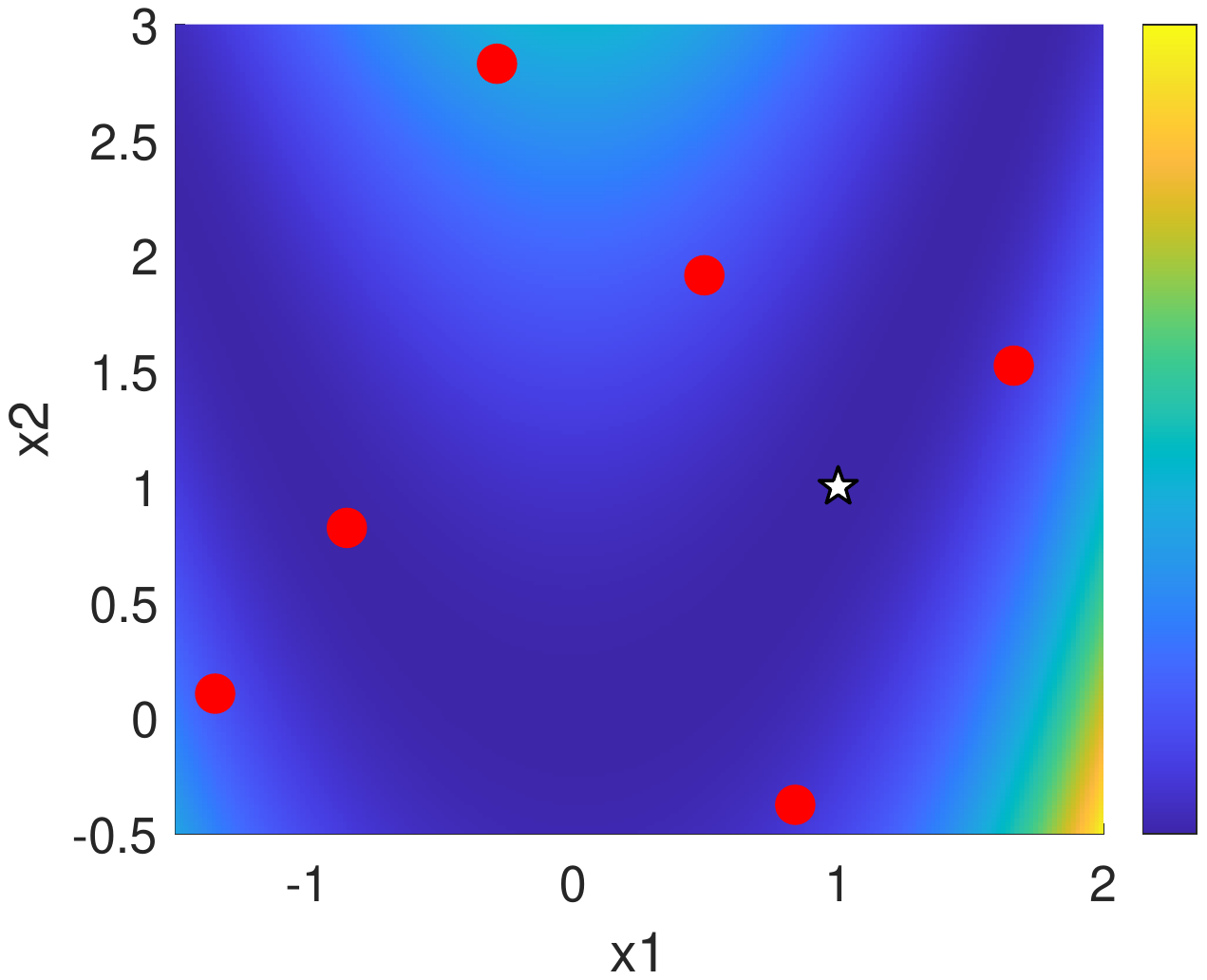}
            \caption{}
        \end{subfigure}
        \hfill
        \begin{subfigure}[b]{0.31\textwidth}
            \centering
%             \includegraphics[scale=.4, width=\textwidth]{./figures/true.pdf}
%            \caption{}
        \end{subfigure}
        \hfill
        \begin{subfigure}[b]{0.31\textwidth}
            \centering
            \includegraphics[scale=.5, width=\textwidth]{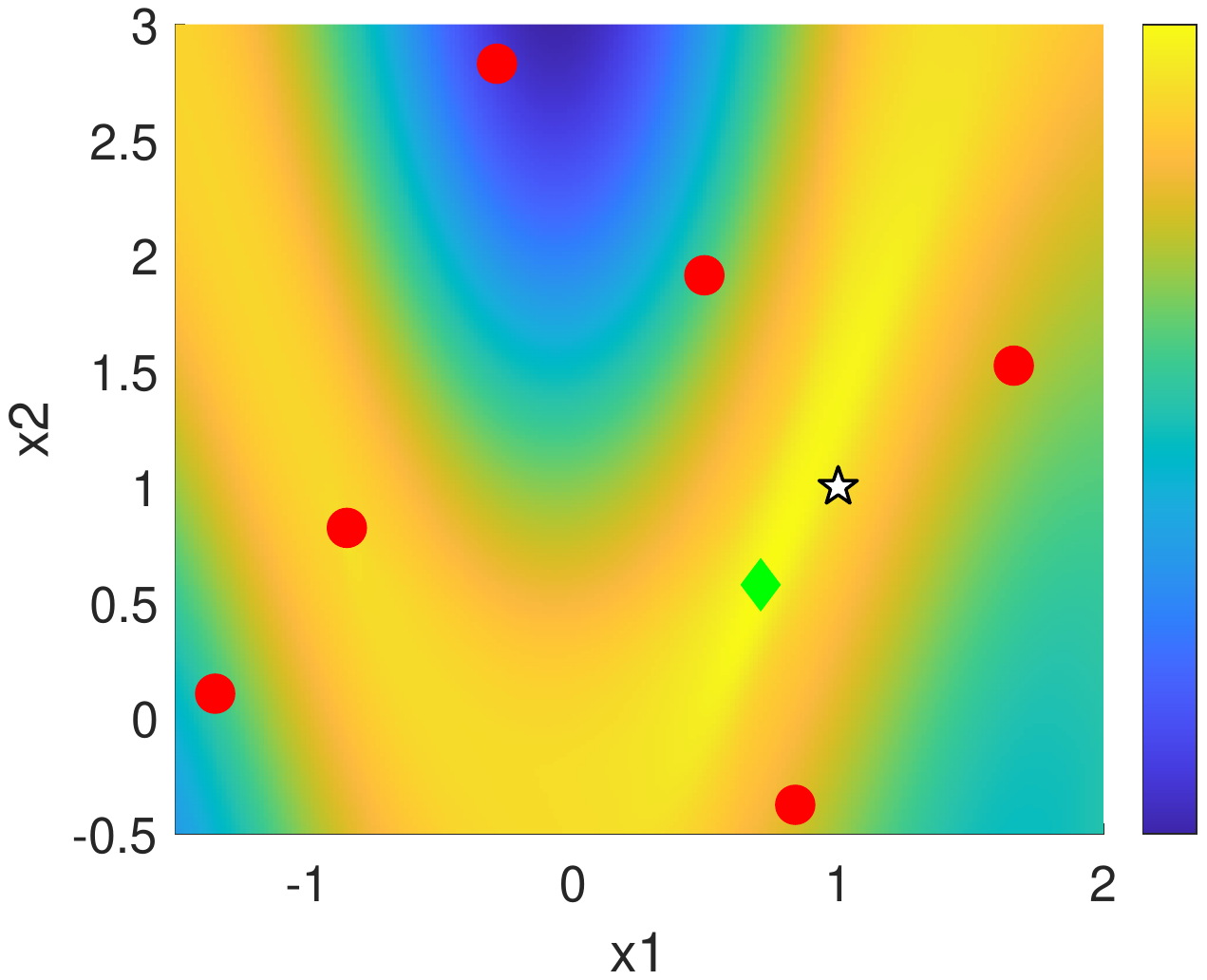}
            \caption{}
        \end{subfigure}
        \caption{Illustrative example comparing the proposed EI-CF and mWB2-CF acquisition functions to the classical EI acquisition function on a modified Rosenbrock problem. The top row shows the (a) mean, (b) variance, and (c) EI for the posterior of the overall objective $\ell(x)$ given observations of $\ell(x_i)$ at six points $x_1, \ldots, x_6$ (shown with red dots in all figures). The second row shows the (d) mean, (e) variance, and (f) EI-CF for the posterior of the composite objective $\ell(x) = f(x,d(x))$ given observations of the black-box function $d(x_i)$ at the same six points. The last row shows the (g) exact shape of the function and (h) mWB2-CF. The white star denotes the true (unknown) global minimizer and the green diamonds denote the maximum of the acquisition functions. }
\label{fig:illustrative_example}
\end{figure}

In Fig. \ref{fig:illustrative_example}, we demonstrate the advantages of exploiting the composite structure of the objective and the proposed mWB2-CF acquisition function on an illustrative example. In particular, we take a grey-box representation of the Rosenbrock function, i.e., $f(x,d(x)) = 100 - x_1^2 + d(x)^2$ where $A = I_2$, $d(x) = x_2-x_1^2$ is the unknown black-box function, and $x = [x_1,x_2]^\top$. The first row of Fig. \ref{fig:illustrative_example} shows the posterior mean (a) and variance (b) of the overall objective $\ell$ for six randomly generated initial samples and the traditional EI acquisition function (c). We see that the predicted variance is low near the evaluated samples (red dots), and this variance grows as $x$ moves away from these points. From Fig. \ref{fig:illustrative_example}c, we see that classical EI is largest near the lowest function value observed (green diamond), which is relatively far away from the true global minimum $f(x^\star, d(x^\star)) = 0$ where $x^\star = [1,1]^\top$ (white star). This is not unexpected since the mean function built for $\ell$ directly does not accurately represent the true surface (g) with these limited set of samples. The second row of Fig. \ref{fig:illustrative_example}, on the other hand, shows the posterior mean (d) and variance (e) of the composite objective $f(x, d(x))$ for the same set of six samples as well as the EI-CF acquisition function (f). It should be noted that the posterior of $f(x, d(x))$ is not normally distributed and thus the mean, variance, and EI-CF evaluations must all be approximated with MC sampling. Not only do we see that the mean function result in a much more accurate representation of true objective, we see significantly less variance in the prediction due to partial knowledge of the function's structure. A direct result of this is that the largest value of EI-CF (green diamond) is much closer to the global minimum, as seen in Fig \ref{fig:illustrative_example}f. However, the max operator in \eqref{eq:ei-cf} results in a flattened EI-CF surface that can be difficult to globally optimize; the mWB2-CF acquisition function (h) clearly addresses this problem by providing useful gradient information throughout the design space while still preserving the global maximum of EI-CF. 

\subsection{Chance constraint-based relaxation of feasible region}

The previous section focused on the choice of acquisition function in the absence of any black- or grey-box constraints ($n_g=0$), which has been the main case of interest in much of the BO literature. As mentioned in the introduction, we take a probabilistic explicit approach to constraint handling that can be generally formulated in terms of chance constraints, which is the focus of this section.

\begin{figure}[ht!]
\centering
\includegraphics[width=0.6\linewidth]{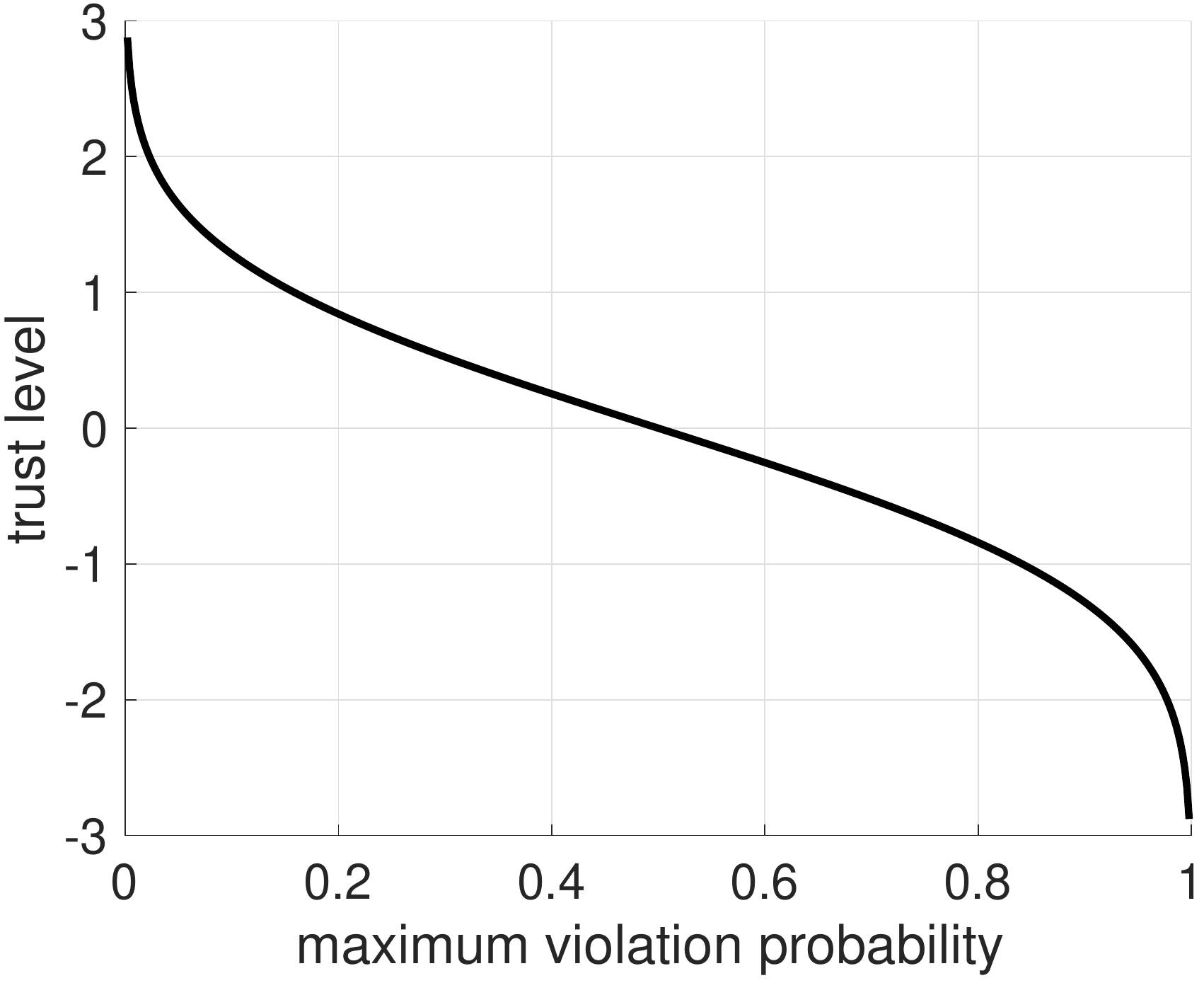}
\caption{The relationship between the maximum violation probability $\epsilon_{i,n}$ and the corresponding trust levels $\tau_{i,n}$ in the probabilistic explicit set representation of the grey-box constraints.}
\label{fig:trust_level_vs_epsilon}
\end{figure}

Given our statistical model of $d$ in  \eqref{eq:dGP}, the constraint function in \eqref{eq:grey-box-problem-g} becomes a multivariate random vector for any particular $x \in \mathcal{X}$. Due to this inherent randomness, we can only enforce these constraints up to a certain probability level. For simplicity, we formulate these as a set of individual chance constraints for each component of $g(\cdot) = [g_1(\cdot), \ldots, g_{n_g}(\cdot)]^\top$ as follows
\begin{align} \label{eq:individual-chance-constraint-g}
\mathbb{P}_n\{ g_i(x, d(Ax)) \leq 0 \} \geq 1-\epsilon_{i,n}, ~~ \forall i \in \{ 1, \ldots, n_g \},
\end{align}
where $\mathbb{P}_n\{ \cdot \}$ is the conditional probability given all available observations $\mathcal{D}_n$ at iteration $n$ and $\epsilon_{i,n} \in (0,1)$ is the maximum allowed probability of constraint violation. 
In similar fashion to EI previously, we see that \eqref{eq:individual-chance-constraint-g} has a closed-form expression in certain special cases. For example, we can derive the following expression when $g_i$ is linearly related to the black-box function, i.e., $g_i(x, d(Ax)) = a_i^\top(x) d(Ax) - b_i(x)$
\begin{align}
a_i^\top(x) \mu_n(Ax) + \Phi^{-1}(1- \epsilon_{i,n}) \| C_n(Ax) a_i \|_2 \leq b_i(x).
\end{align}
Unfortunately, such closed-form expressions are not available in the general nonlinear case so we again need to develop an efficient approximation strategy. One approach would be to again try a MC sampling-based approach by recognizing that 
\begin{align} \label{eq:chance-constraint-saa}
\mathbb{P}_n\{ g_i(x, d(Ax)) \leq 0 \} &= \mathbb{E}_n \{ \mathbb{I}\{ g_i(x, d(Ax)) \leq 0 \} \}, \\\notag
&= \mathbb{E}_\xi \{ \mathbb{I}\{ g_i(x, \mu_n(Ax) + C_n(Ax) \xi) \leq 0 \} \}, \\\notag
&\approx \frac{1}{M}\sum_{i=1}^M \mathbb{I}\{ g_i(x, \mu_n(Ax) + C_n(Ax) \xi^{(i)}) \leq 0 \},
\end{align}
where $\mathbb{I}\{ L \}$ is the indicator function that is equal to 1 if the logical proposition $L$ is true and 0 otherwise and $\{ \xi^{(i)} \}_{i=1}^M$ is a set of $M$ i.i.d. samples of $\xi$. Since the indicator function in \eqref{eq:chance-constraint-saa} is non-smooth, it substantially increases the complexity of the enrichment sub-problem due to the inclusion of binary variables to enforce some fraction of the constraints to hold \cite{pagnoncelli09}. A much simpler approach is to rely on a moment-based approximation of \eqref{eq:individual-chance-constraint-g} that can be derived from a first-order Taylor series expansion of $g_i(x,y)$ in its second argument around $\hat{y} = \mu_n(Ax)$
\begin{align}
g_i(x, y) \approx g_i(x, \hat{y}) + \nabla_y g_i(x, \hat{y})(y - \hat{y}).
\end{align} 
By substituting $y = d(Ax)$ in the expression above, we see that the posterior distribution of $g_i(x, d(Ax))$ is approximately normally distributed for any fixed $x \in \mathcal{X}$
\begin{align}
g_i(x, d(Ax)) | \mathcal{D}_n \approx \mathcal{N}( \mu_{g_i, n}(x) , \sigma^2_{g_{i}, n}(x) ),
\end{align}
%g_i(x, d(Ax)) \approx \mathcal{N}( g_i(x, \mu_n(Ax)), \nabla_y g_i(x, \mu_n(Ax)) \Sigma_n(Ax) \nabla_y g_i(x, \mu_n(Ax))^\top ),
where 
\begin{subequations} \label{eq:musigmag}
\begin{align}
\mu_{g_i, n}(x) &= g_i(x, \mu_n(Ax)), \\
\sigma^2_{g_{i}, n}(x) &= \nabla_y g_i(x, \mu_n(Ax)) \Sigma_n(Ax) \nabla_y g_i(x, \mu_n(Ax))^\top,
\end{align}
\end{subequations}
are the approximate mean and variance for the probabilistically modeled constraint function $g_i$ at iteration $n$, respectively. Under this approximation, the chance constraints \eqref{eq:individual-chance-constraint-g} can be substantially simplified, which enables us to define the probabilistic explicit set $\mathcal{X}_n$ as follows
\begin{align} \label{eq:Xnset}
\mathcal{X}_n = \{ x \in \mathcal{X} : \mu_{g_i, n}(x) + \tau_{i,n} \sigma_{g_{i}, n}(x) \leq 0, ~~ \forall i \in \{1, \ldots, n_g \} \}, 
\end{align}
where $\tau_{i,n} = \Phi^{-1}(1- \epsilon_{i,n})$ can be interpreted as ``trust'' levels in the constraints at the current iteration number $n$. We see $\tau_{i,n}$ is explicitly related to the allowed violation probability $\epsilon_{i,n}$, with the shape of this curve being depicted in Fig. \ref{fig:trust_level_vs_epsilon}. For $\epsilon_{i,n}$ values less than 50\%, we see that $\tau_{i,n}$ is positive such that the variance term $\tau_{i,n} \sigma_{g_{i}, n}(x)$ effectively ``backs off'' of the mean prediction of the constraint. This is an overall restriction of a simple deterministic prediction of the feasible region, i.e., $g(x, \mu_n(Ax)) \leq 0$, which is recovered in the case that $\epsilon_{i,n} = 0.5$. For $\epsilon_{i,n}$ values greater than 50\%, on the other hand, we see that $\tau_{i,n}$ becomes negative, so that \eqref{eq:Xnset} is a relaxation of the nominal prediction of the feasible region; the basic idea of this relaxation is to try to ensure that the predicted feasible region contains the true one with high probability and is similar to the notion of the upper trust bounds (UTBs) introduced in \cite{priem19}. In fact, it can be shown that \eqref{eq:Xnset} reduces to the UTB approach in the case of fully black-box constraint functions. 

We can think of $\epsilon_{i,n}$ (or equivalently $\tau_{i,n}$) as tuning parameters of the algorithm. In fact, one can imagine developing various adaptive tuning mechanisms to update $\tau_{i,n}$ at each iteration $n$ to control the level of exploration in the predicted feasible region. In the early stages of Algorithm \ref{alg:gp-sampling}, we expect the GP to provide a poor prediction of the feasible region such that larger values of $\tau_{i,n}$ would be preferred to allow for more exploration instead of overly sampling near a feasible local minimum. As $n$ increases, however, we expect the GP model to become more accurate and, as our confidence increases, we should shrink $\tau_{i,n}$ toward zero to accelerate convergence. Using this rule-of-thumb, we found that the following linear update rule worked well in our numerical tests
\begin{align} \label{eq:tau-tuning}
\tau_{i,n} = -3\left(1 - \frac{n}{N} \right),
\end{align}
which starts from an initial condition of $-3$ (corresponds to a high confidence that the true region is contained in $\mathcal{X}_n$) and converges to 0 (corresponds to the nominal GP prediction) at the final iteration. A variety of other update strategies for the trust levels could be considered including ones that exhibit the opposite behavior or even further increase $\tau_{i,n}$, which we plan to explore more in future work. 

Since $\mathcal{X}_n$ in \eqref{eq:Xnset} is implicitly defined by a collection of nonlinear constraints, the stochastic gradient ascent algorithm in \eqref{eq:sga} becomes significantly more expensive as a direct result of the projector operator $\Pi_{\mathcal{X}_n}$ that requires a non-convex optimization to be solved at every iteration. A more practical alternative is to resort to a sample average approximation (SAA) \cite{kleywegt02} of the enrichment sub-problem \eqref{eq:optimize-acqusition} that reduces the stochastic optimization problem to a deterministic one. Using the proposed $\text{mWB2-CF}_n(x)$ acquisition function in \eqref{eq:modified-WB2CF} and probabilistic explicit set representation in \eqref{eq:Xnset}, the SAA approximation of \eqref{eq:optimize-acqusition} can be formulated as
\begin{subequations} \label{eq:saa}
\begin{align}
x_{n+1} = \argmax_{x} &~~ \frac{1}{M} \sum_{i=1}^M \left( s_n [\ell_n^\star - \ell_n^{(i)}(x)]^+ - \ell_n^{(i)}(x) \right), \\\label{eq:saa-constraints}
\text{s.t.} &~~ \mu_{g_i, n}(x) + \tau_{i,n} \sigma_{g_{i}, n}(x) \leq 0, ~~ \forall i \in \{ 1, \ldots, n_g \}, \\
&~~ x^L \leq x \leq x^U, 
\end{align}
\end{subequations}
where $\ell_n^{(i)}(x) = f(x, \mu_n(Ax) + C_n(Ax) \xi^{(i)})$ for all $i = 1,\ldots, M$; $\xi^{(1)}, \ldots, \xi^{(M)}$ is a set of $M$ i.i.d. realizations of the random variable $\xi \sim \mathcal{N}(0, I_{n_y})$; and the incumbent definition must be modified to account for constraints as follows
\begin{align} \label{eq:constrained-incumbent}
\ell_n^\star = \begin{cases}
\infty,  &\text{if no points are feasible}, \\
\displaystyle\min_{i \in \{ 1,\ldots, n \}} f(x_i, y_i)~~\text{s.t.} ~~ g(x_i,y_i) \leq 0, &\text{otherwise}.
\end{cases}
\end{align}
In practice, when the incumbent is infinite implying no feasible points exist, we set $s_n=0$ such that the objective reduces to a mean-based acquisition function. 

The key distinction of SAA (relative to the stochastic gradient ascent method) is that the \textit{same} set of realizations are used for the different iterates of $x$ computed during the optimization process such that \eqref{eq:saa} is fully deterministic. This implies that we can use any of the state-of-the-art (derivative-based) nonlinear programming (NLP) solvers to efficiently locate locally optimal solutions to this problem. Due to the fast local convergence rate of derivative-based NLP solvers, such as IPOPT \cite{biegler09}, we can easily utilize them in a multi-start optimization procedure as a heuristic to finding the global maximizer of \eqref{eq:saa}. It is also important to note that this heuristic is not a major limitation, as even approximate solutions to \eqref{eq:saa} still produce valuable samples that are more likely to result in improvement of the current incumbent relative to alternatives such as random or space filling designs. 

\begin{remark} \label{rem:2}
There has been a substantial amount of theory established for the SAA method applied to stochastic optimization problems (see, e.g., \cite{kim15} for detailed discussions on these results). Under certain assumptions on the objective and constraint functions, the optimal solution and objective estimate obtained from \eqref{eq:saa} converge to their respective true values \textit{in distribution} at a rate of $O(1/\sqrt{M})$. Furthermore, various procedures have been developed to estimate the optimality gap between the solution of \eqref{eq:saa} and the true solution. This optimliaty gap estimator and its corresponding variance estimated from the MC standard error formula can be used to decide if the approximate solution is trustworthy or if the number of samples $M$ needs to be increased, as discussed in \cite[Section 3]{kleywegt02}.
\end{remark}

\begin{remark} \label{rem:3}
Although we focused on inequality constraints \eqref{eq:grey-box-problem-g} for simplicity, we can handle nonlinear equality constraints of the form $h(x,y) = 0$ similarly to \eqref{eq:Xnset}. In particular, we can transform $h(x,y) = 0$ into two-sided inequalities $h(x, y) \leq 0$ and $h(x, y) \geq 0$ to arrive at
\begin{align} \label{eq:equality-constraints}
| \mu_{h_i, n}(x) | \leq -\tau_{i,n} \sigma_{h_{i}, n}(x), ~~ \forall i \in \{ 1, \ldots, n_h \},
\end{align}
where $\mu_{h_i, n}(x)$ and $\sigma_{h_{i}, n}(x)$ are the approximate mean and standard deviation for the $i^\text{th}$ component of $h$ at iteration $n$, which are defined similarly to \eqref{eq:musigmag}. Note that the constraints \eqref{eq:equality-constraints} become infeasible for any positive $\tau_{i,n}$, implying that update rules that keep these values negative, such as \eqref{eq:tau-tuning}, must be used. Furthermore, when $\tau_{i,n}=0$, \eqref{eq:equality-constraints} can be simplified to $h(x, \mu_n(Ax)) = 0$, which corresponds to the best nominal prediction of the grey-box equality constraints. 
\end{remark}

\begin{remark} \label{rem:4}
We focused on an individual chance constraint formulation \eqref{eq:individual-chance-constraint-g}, as opposed to the following joint chance constraint representation
\begin{align} \label{eq:jcc}
\mathbb{P}_n\{ g(x, d(Ax)) \leq 0 \}  = \mathbb{P}_n\{ g_i(x, d(Ax)) \leq 0, \forall i \in \{ 1, \ldots, n_g\} \} \geq 1 - \epsilon,
\end{align}
for two reasons. First, individual chance constraints can be dealt with in a more computationally efficient manner since they do not require the construction of the full multivariate distribution of $g(x, d(Ax)) | \mathcal{D}_n$. Second, the choices of $\{ \epsilon_{i,n} \}$ (that can be equivalently mapped to the trust levels $\{ \tau_{i,n} \}$) remain tuning parameters in COBALT. Although \eqref{eq:jcc} involves only a single tuning parameter $\epsilon$, it is not immediately obvious how to select this value to achieve good performance in practice. Thus, we opt for the approach in \eqref{eq:tau-tuning} due to the ease of interpretation of the impact on the selected $x_{n+1}$ to changes in the trust levels. We plan to study the impact of the individual versus joint chance constraint representation more in our future work. It is interesting to note, however, that the collection of individual violation probabilities $\{ \epsilon_{i,n} \}$ can be used to derive a conservative estimate on the joint satisfaction probability using Boole's inequality \cite{paulson20b}, i.e., $\epsilon \leq \sum_{i=1}^{n_g} \epsilon_{i,n}$, though this bound is known to be conservative (unless the events are nearly independent). 
\end{remark}

%%%%%%%%%%%%%%%%%%%%%%%%%
\section{Numerical Examples}
\label{sec:numerical-example}

Several numerical experiments were carried out to test the performance of the COBALT algorithm under a variety of different conditions. In all cases, we assume that noiseless objective and constraint functions are available, however, the algorithm is capable of identifying the noise term in the GP, as discussed in Section \ref{subsec:GPsingle}, and so is directly applicable to cases with noisy observations. 

\subsection{Implementation of the COBALT algorithm}

We now outline our implementation of the proposed COBALT algorithm that was used to generate results in this section. The code, which is implemented in Matlab, is freely available online at \cite{paulson21}. The first major step is construction of the GP regression models for the black-box functions, which require the hyperparameters of the covariance function to be estimated, as discussed in Section \ref{subsec:GPsingle}. This is achieved by solving the MLE problem in \eqref{eq:MLE} using the DIRECT search algorithm from \cite{finkel03}. The best hyperparameter value found by the initial DIRECT search is then used to initialize Matlab's local \texttt{fmincon} solver. The other main GP calculations, such as evaluating the mean and variance functions in \eqref{eq:mean-variance-s}, were performed with the Gaussian Process for Machine Learning (GPML) toolbox developed by Rasmussen and Nickisch \cite{rasmussen10}. The GPML toolbox implements the Mat\'ern and squared exponential covariance functions in \eqref{eq:Matern-cov} and \eqref{eq:SE-cov}, respectively. We selected the squared exponential kernel by default in all of our tests, though this can easily be modified in the code. The SAA-based enrichment sub-problem \eqref{eq:saa}, with $M=100$ as a default, was solved using IPOPT \cite{biegler09}, with the CasADi \cite{andersson19} automatic differentiation toolbox being used to supply exact first- and second-order derivatives to IPOPT. 

Note that the current version of COBALT \cite{paulson21} is constructed to be modular, as it implements a \texttt{GreyboxModel} class object that has several helper methods that can be useful for algorithm testing and comparisons. With the exception of the final case study that requires specialized simulation code to run, all other test problems have been included in the initial code release. 

\subsection{Optimization test problems and performance assessment method}

We test the performance of COBALT on a diverse set of optimization test problems commonly used for benchmarking global optimization algorithms, which we slightly modify to fit the proposed grey-box structure in this paper. A summary of the set of seven test problems is provided in Table \ref{tab:testproblems}. The exact equations for the test problems in the form of \eqref{eq:grey-box-problem}, along with their global solutions, are provided in Appendix \ref{appendix:A}. The first three problems involve highly nonlinear (and some multi-modal) composite objective functions of varying dimension with only box constraints on the decision variables. The next three problems all involve nonlinear constraints, with at least one or more composite grey-box constraint functions. While the Toy-Hydrology problem has a known objective and one grey-box constraint function, the Rosen-Suzuki and Colville problems have a mixture of grey-box objective and constraint functions. The last problem is a realistic engineering problem that is related to parameter estimation in complex bioreactor models. 

\begin{table}[ht!]
\caption{Overview of the characteristics of the collection of seven optimization test problems considered in this work along with their corresponding source. The detailed formulation of each test problem is provided in Appendix \ref{appendix:A}.}
\begin{center}
\begin{tabular}{l c c c c c l} \hline
Name &$n_x$ &$n_y$ &$n_z$ &$n_g$ &Equation \# &Reference \\
\hline
Goldstein-Price &2 &2 &2 &0 &\eqref{eq:GSP} &Dixon, \cite{dixon78} \\
Rastrigin &3 &1 &1 &0 &\eqref{eq:Rastrigin} &Rastrigin, \cite{rastrigin74} \\
Rosenbrock &6 &4 &4 &0 &\eqref{eq:Rosenbrock} &Rosenbrock, \cite{rosenbrock60} \\
Toy-Hydrology &2 &1 &1 &2 &\eqref{eq:ToyHydrology} &Gramacy et al., \cite{gramacy16} \\
Rosen-Suzuki &4 &2 &2 &3 &\eqref{eq:RosenSuzuki} &Hock et al., \cite{hock80} \\
Colville &5 &4 &4 &6 &\eqref{eq:Colville} &Rijckaert et al., \cite{rijckaert78} \\
DFBA-MLE &6 &24 &6 &0 &\eqref{eq:dfba-substitutions} &Paulson et al., \cite{paulson19} \\
\hline
\end{tabular}
\label{tab:testproblems}
\end{center}
\end{table}

We compare the performance of COBALT (mWB2-CF) with four other acquisition functions: expected improvement (EI), probability of improvement (PI), the expected improvement for composite functions (EI-CF), and random search (Random). The EI and PI acquisitions were run using the \texttt{bayesopt} function from the Statistics and Machine Learning Toolbox in Matlab. The EI-CF acquisition was implemented in a similar fashion to \eqref{eq:saa}, with the mWB2-CF objective function from \eqref{eq:modified-WB2CF} replaced by \eqref{eq:ei-cf}. For all problems and methods, an initial set of evaluations is performed using $\max\{ 3, n_z+1 \}$ points chosen with Latin hypercube sampling (LHS). We use the base-10 logarithm of the best-sample simple regret as our performance metric, which is defined as follows
\begin{align} \label{eq:regret}
\text{Log10-Regret}_n = \log_{10}(\ell^\star_n - \ell^\star_\text{true}), 
\end{align}
where $\ell^\star_\text{true}$ is the global minimum of the exact optimization problem (see Appendix \ref{appendix:A}). The regret is a measure of how far off the best currently feasible sample, defined in \eqref{eq:constrained-incumbent}, is from the true global optimum at each iteration of the algorithm. The logarithm is incorporated to account for the largely different scales that can occur depending on the specifics of the objective function and the initialization points. Since $\text{Log10-Regret}_n$ depends on the randomly selected initial points, showing results for a single initialization is not very informative. As an alternative, we repeat every experiment 50 times to estimate the average $\text{Log10-Regret}_n$ for each algorithm. Error bars are computed by estimating the confidence intervals as 1.96 times the standard deviation divided by the square root of the number of repeats. Furthermore, a detailed overview of the computational statistics of COBALT (relative to traditional constrained BO) is provided in Appendix \ref{appendix:B}.

\subsection{Results for box-constrained global optimization test problems}

\begin{figure}[tb!]
\centering
\includegraphics[width=0.6\linewidth]{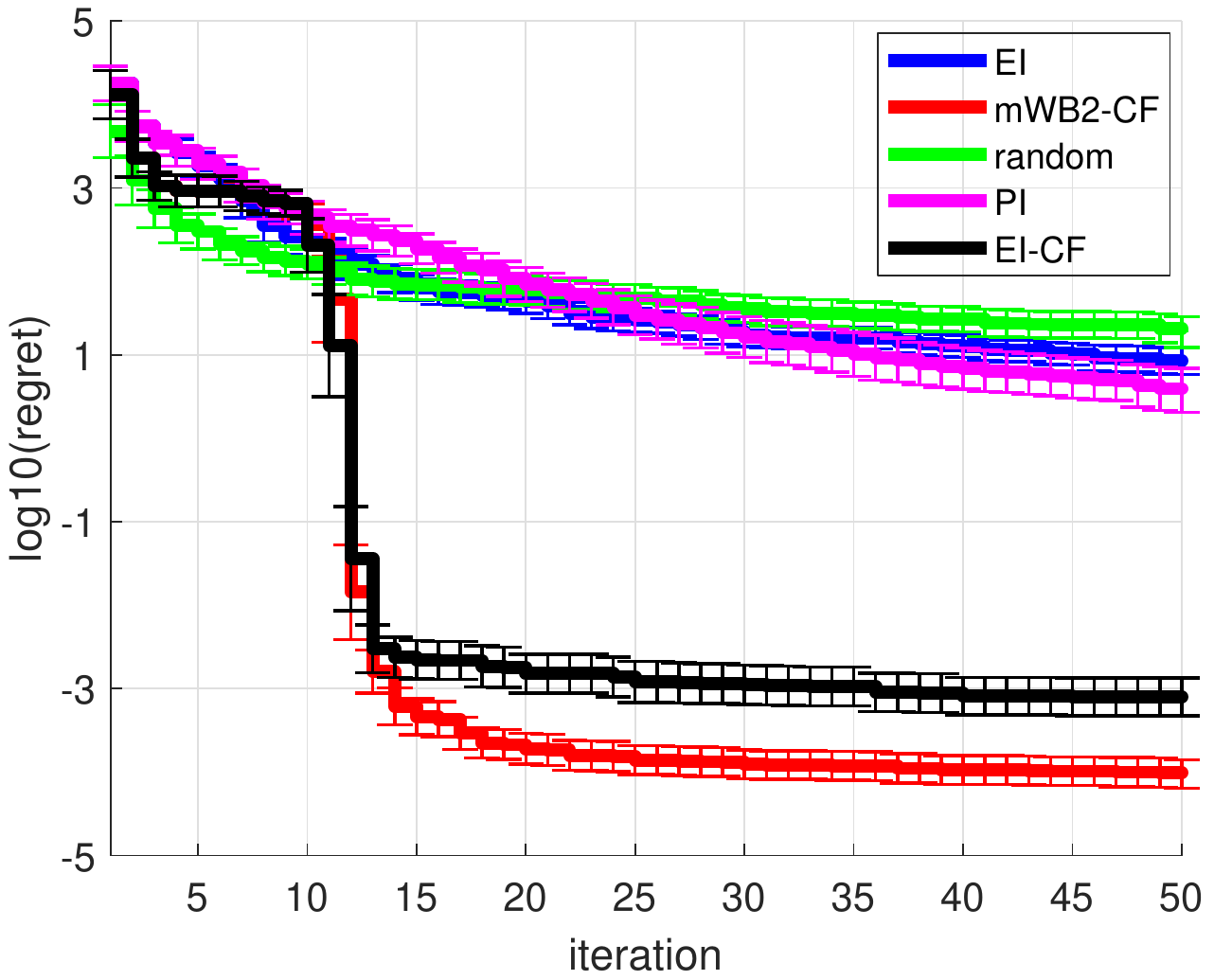}
\caption{Expected $\log10(\text{regret})$ for the 2d Goldstein-Price test function, with approximate confidence region shown via error bars, estimated from 50 independent realizations.}
\label{fig:GSP}
\end{figure}

\begin{figure}[tb!]
\centering
\includegraphics[width=0.6\linewidth]{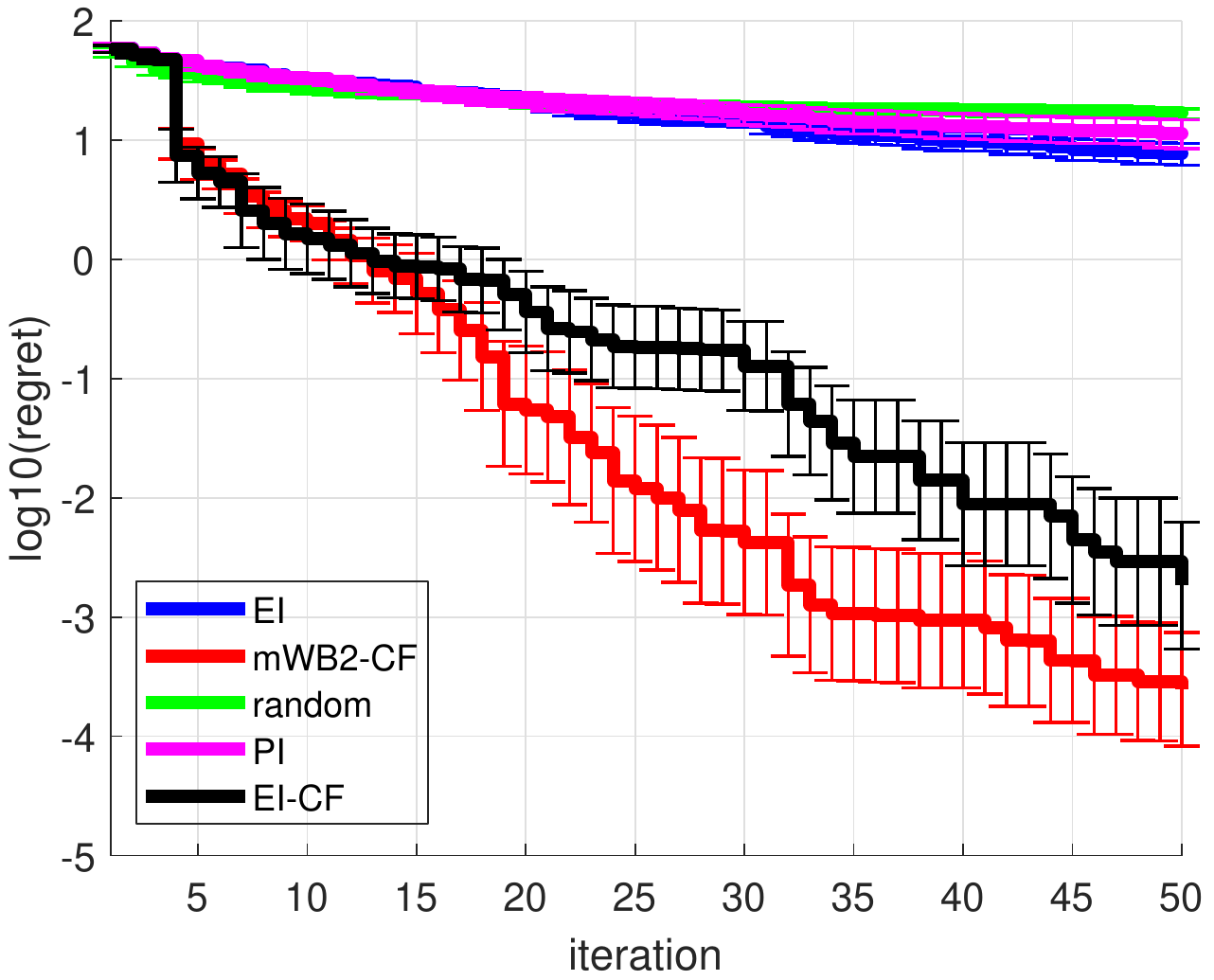}
\caption{Expected $\log10(\text{regret})$ for the 3d Rastrigin test function, with approximate confidence region shown via error bars, estimated from 50 independent realizations.}
\label{fig:Ras3}
\end{figure}

\begin{figure}[tb!]
\centering
\includegraphics[width=0.6\linewidth]{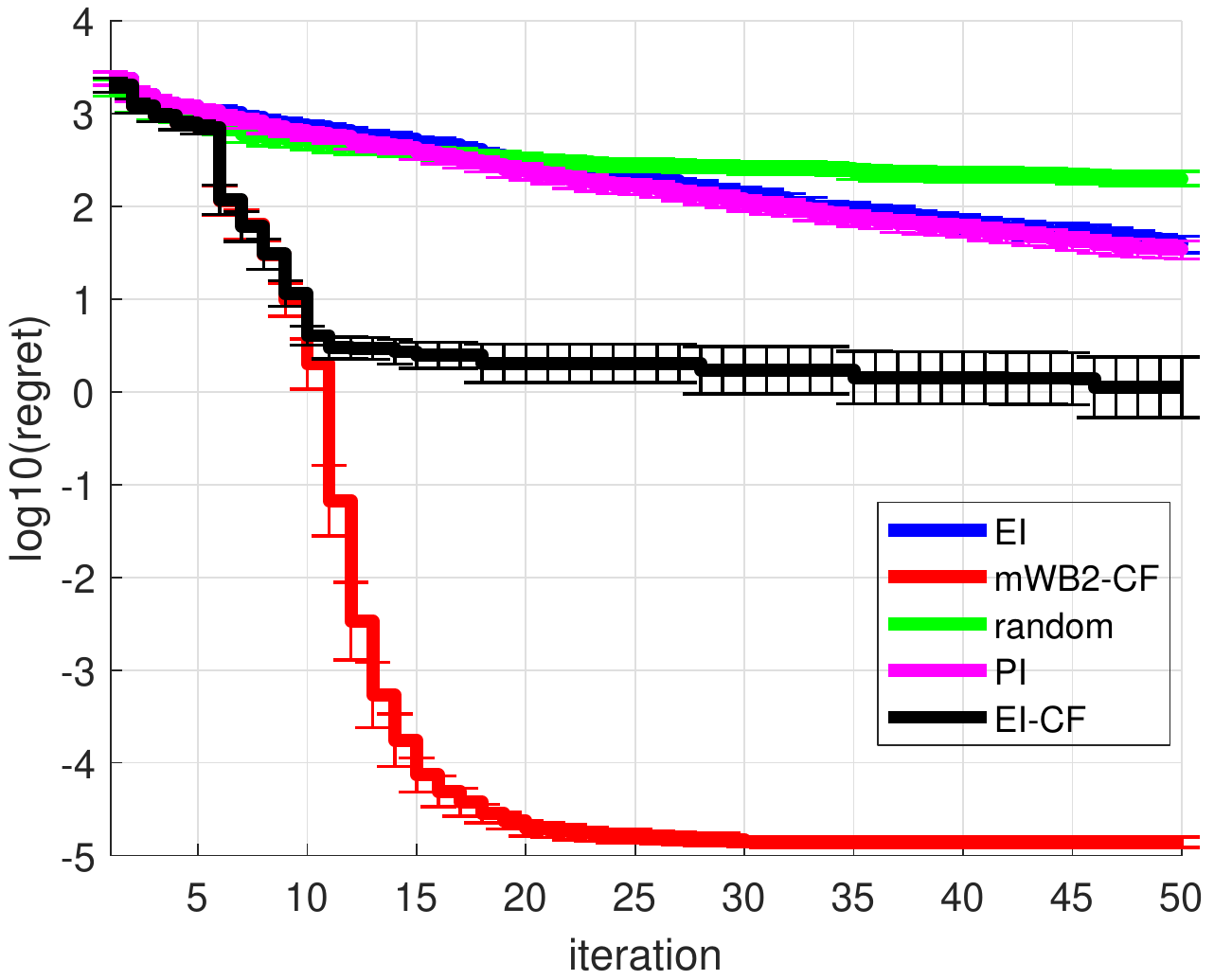}
\caption{Expected $\log10(\text{regret})$ for the 6d Rosenbrock test function, with approximate confidence region shown via error bars, estimated from 50 independent realizations.}
\label{fig:R6}
\end{figure}

We first discuss the results for the set of box-constrained test optimization problems. Figs. \ref{fig:GSP}, \ref{fig:Ras3}, and \ref{fig:R6}, respectively, show the expected $\text{Log10-Regret}_n$ over 50 replications for the Goldstein-Price, Rastrigin, and Rosenbrock functions. We clearly see that the COBALT algorithm (mWBS-CF) outperforms all other tested acquisition functions by up to 3 orders of magnitude. It is interesting to note that, even though the Rastrigin function in particular has a large number of local minima, both grey-box acquisition functions do not get ``stuck'' and are able to make much faster progress than their fully black-box counterparts. For the smaller dimensional problems (Goldstein-Price and Rastrigin), we see that EI-CF performs only slightly worse than mWB2-CF; however, for the six dimensional Rosenbrock problem, mWB2-CF results in significant improvement over EI-CF. This suggests that the proposed choice of acquisition function (and its sub-optimization routine) has a large influence on performance, and that the gradient of mWB2-CF plays an important role in selecting samples that are more likely to be optimal. 
We also observe that the confidence intervals for the average regret are higher for the Rastrigin problem compared to the Goldstein-Price or Rosenbrock problems. First, it should be noted that regret is plotted on a logarithmic scale so that the actual variance is small for all of these problems. Second, the Rastrigin objective function is quite ``bumpy'' due to the inclusion of a cosine term that is difficult to resolve with a small number of samples. As such, the performance of any BO algorithm will more strongly depend on the initial set of randomly selected samples that naturally leads to more run-to-run variability in the tested algorithms. 

\subsection{Results for nonlinearly constrained global optimization test problems}
%\subsection{Performance in the absence of a feasible starting point}

\begin{figure}[tb!]
\centering
\includegraphics[width=0.6\linewidth]{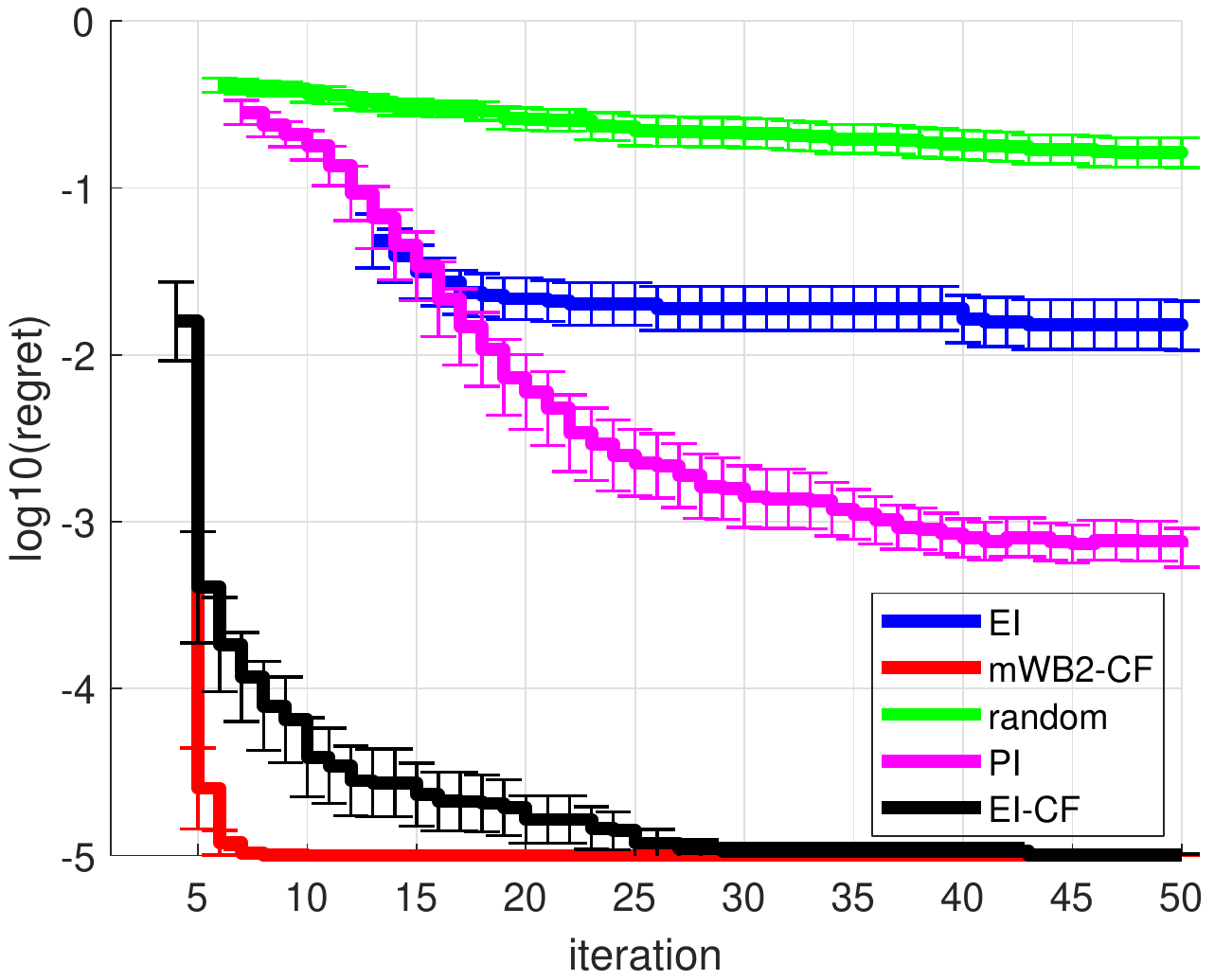}
\caption{Expected $\log10(\text{regret})$ for the 2d Toy-Hydrology test function, with approximate confidence region shown via error bars, estimated from 50 independent realizations.}
\label{fig:ToyH}
\end{figure}

\begin{figure}[tb!]
\centering
\includegraphics[width=0.6\linewidth]{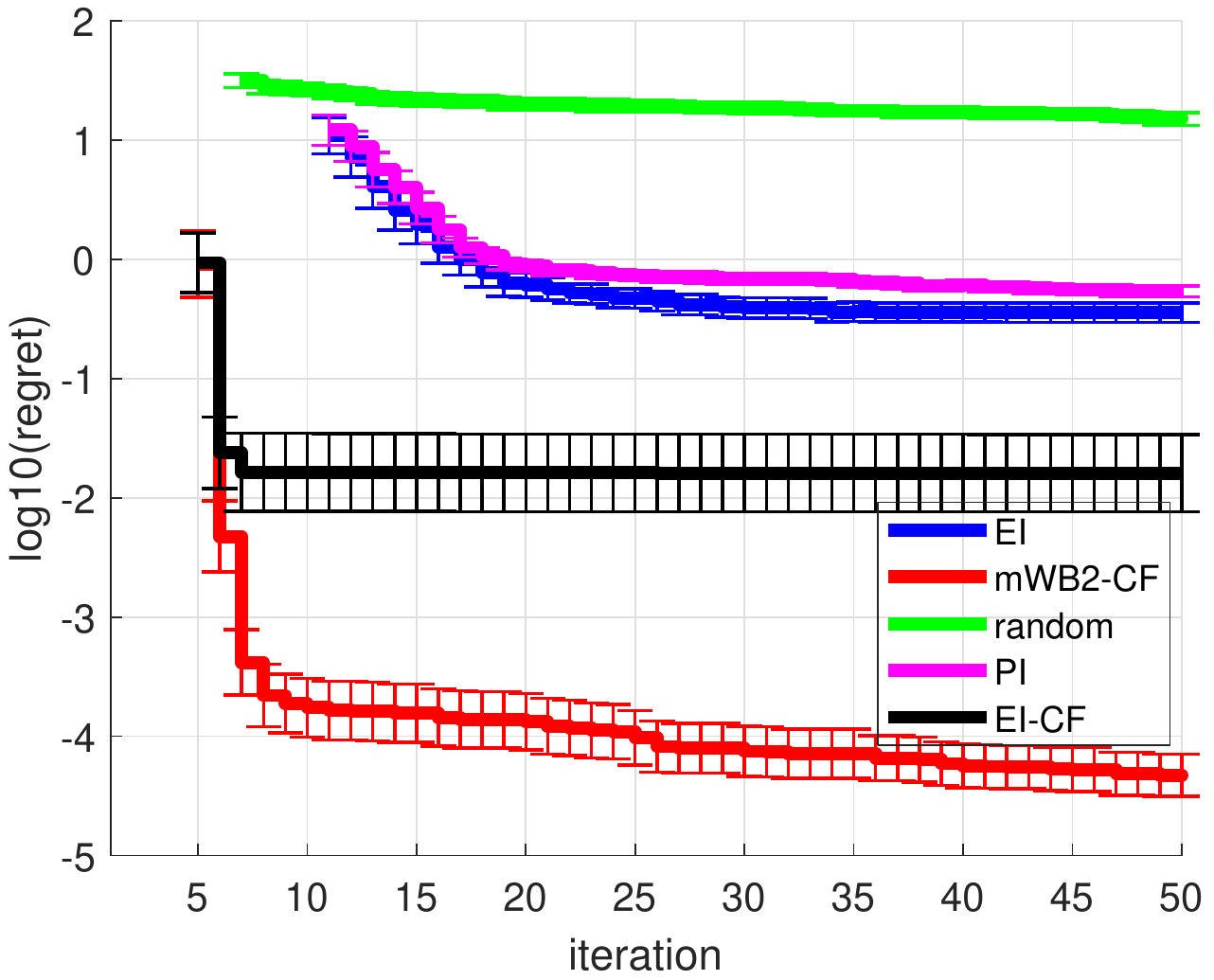}
\caption{Expected $\log10(\text{regret})$ for the 4d Rosen-Suzuki test function, with approximate confidence region shown via error bars, estimated from 50 independent realizations.}
\label{fig:RosenSu}
\end{figure}

\begin{figure}[tb!]
\centering
\includegraphics[width=0.6\linewidth]{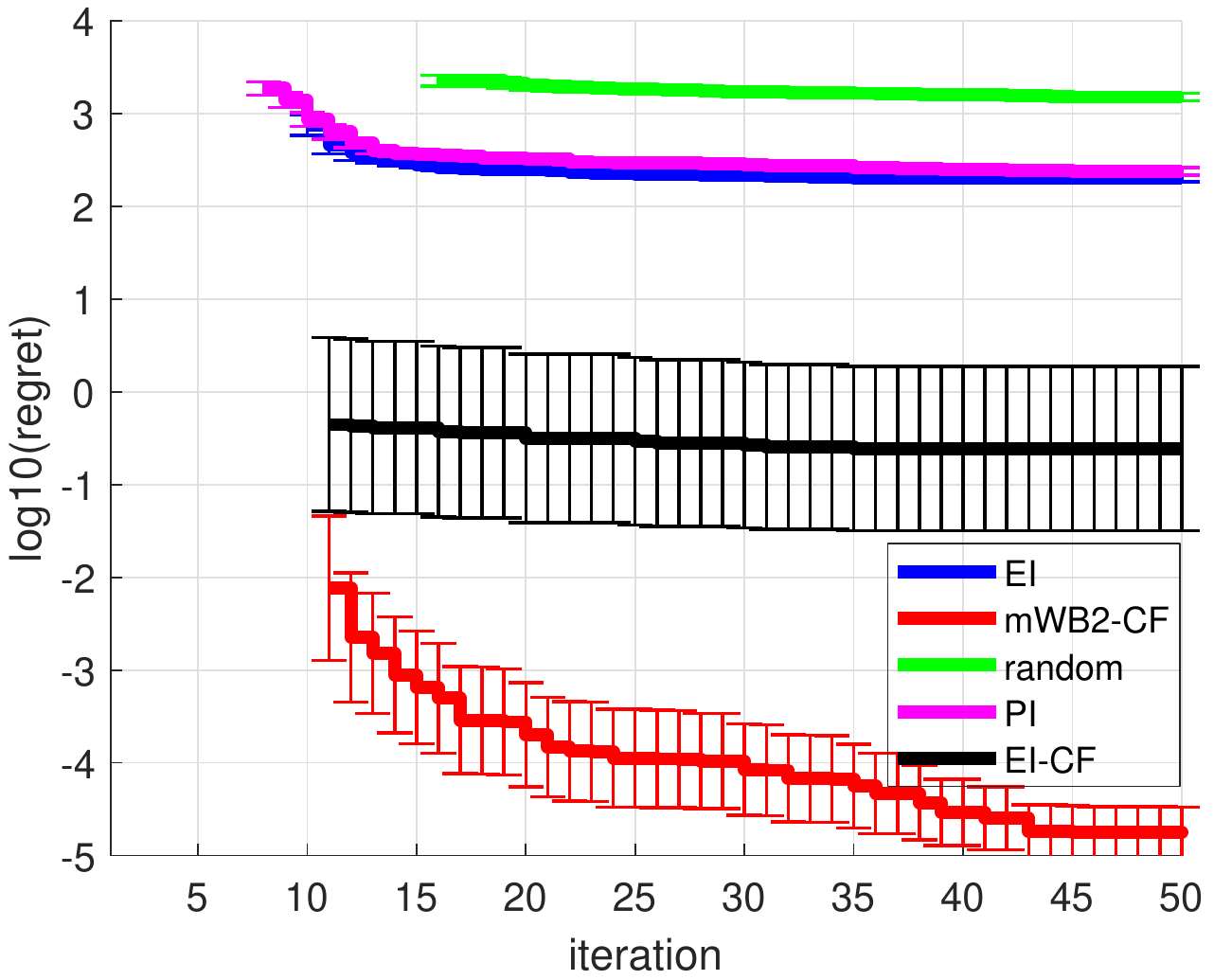}
\caption{Expected $\log10(\text{regret})$ for the 5d Colville test function, with approximate confidence region shown via error bars, estimated from 50 independent realizations.}
\label{fig:Colville}
\end{figure}

Now, we highlight the performance of the COBALT algorithm (mWBS-CF) on the test optimization problems with highly nonlinear (grey-box) composite constraint functions. In particular, the expected value of $\text{Log10-Regret}_n$ over 50 replications of the randomly selected initialization points for the Toy-Hydrology, Rosen-Suzuki, and Colville problems are shown in Figs. \ref{fig:ToyH}, \ref{fig:RosenSu}, and \ref{fig:Colville}, respectively. Recall that COBALT uses a probabilistic explicit constraint handling method, as shown in \eqref{eq:saa-constraints}, whereas \texttt{bayesopt} uses an implicit approach that effectively multiplies the unconstrained acquisition function by the probability of constraint sanctification (see, e.g., \cite{gardner14} for additional details). 

Since the initial samples are selected randomly, there is no guarantee that Algorithm \ref{alg:gp-sampling} (in the grey- or black-box case) will start from a feasible point. As the algorithm progresses and a better model of the black-box components of the model are obtained, we expect that it becomes more likely that a feasible point is identified. To clearly represent this behavior in Figs. \ref{fig:ToyH}--\ref{fig:Colville}, we do not plot the $\text{Log10-Regret}_n$ values for any $n$ in which an infeasible incumbent value exists in any of the 50 replicate runs. As such, the starting iteration (e.g,. $n=5$ for mWB2-CF in Fig. \ref{fig:ToyH}) on these plots is an estimate of the worst-case number of iterations required to find a point within the unknown feasible domain. For the lower dimensional problems (Toy-Hydrology and Rosen-Suzuki), we see that COBALT is not only able to find feasible points faster than \texttt{bayesopt}, but shows a significantly faster decrease in the regret \eqref{eq:regret} to near-zero values. This is likely due to the fact that the black-box components of the grey-box model are defined in terms of fewer variables $n_z/n_x=0.5$, so that a more accurate representation of the feasible region can be constructed with less data. The opposite trend, however, is observed for the higher dimensional problem (Colville), which is likely due to two factors: (i) a larger number of variables interacting with the black-box model $n_z/n_x = 0.8$ and (ii) the relaxation of predicted feasible domain due to the negative-valued trust levels, which gradually shrink as iterations increase \eqref{eq:tau-tuning}. The increased exploration achieved due to this relaxation does provide a substantial gain in the regret, which is seven order of magnitudes better than \texttt{bayesopt} for the mWB2-CF acquisition function. 
Similarly to the behavior previously observed for the Rastrigin problem, we see that the Colville problem results in higher variance in the regret than the Toy-Hydrology and Rosen-Suzuki problems. In this case, the increased regret variance can be attributed to a higher dimensional $d(z)$ with strong nonlinear interaction terms that are harder to estimate. Furthermore, the third component of the unknown function $d_3(z)$ plays an important role in one of the active constraints $g_5(x,y)$ (see Appendix \ref{subappendix:colville}). Compare this to the Rosen-Suzuki problem, for example, which only has two unknown quadratic functions and one of these functions appears in an inactive constraint (so the solution is overall less sensitive to the estimated GP model for $d(z)$).
%{\color{blue}
%Looking at the confidence intervals for the average regret, we see that the Rastrigin and Colville problems show higher variance than what was observed for the other test problems. First, we note that the average regret is plotted on a log scale so that the actual value of the variance for these problems is small. Second, both the Rastrigin and Colville problems have a relatively large number of local solutions, meaning it is easy for the algorithms to over-exploit for some number of iterations, especially when a limited amount of function evaluations are available. In such cases, we expect the regret to show a stronger dependence on the initial set of randomly selected samples, implying the variance of the regret decays more slowly with the number of iterations. 
%}

Lastly, we note that, even though EI-CF performs similarly to mWB2-CF in the Toy-Hydrology problem, we see significant advantages of mWB2-CF in the Rosen-Suzuki and Colville problems. This can be attributed to the fact that EI-CF has zero gradient (see \eqref{eq:ei-cf-gradient} and \eqref{eq:gamma-n}) for a significant portion of the feasible space once a reasonably good solution has been found. This effectively prevents EI-CF from making good progress during the later iterations (due to the difficulty in finding the true global maximizer) that is clearly mitigated by the switch to mWB2-CF. 

\subsection{Results for Bayesian calibration of genome-scale bioreactor model}

Our final case study is focused on a realistic engineering problem related to the calibration of expensive computer model to data. In particular, we are interested in estimating parameters in the substrate uptake kinetics of a genome-scale bioreactor model from batch measurements of concentrations of extracellular metabolites and biomass. This problem was originally proposed in \cite{paulson19} and solved using a surrogate-based optimization approach that required 1200 expensive model evaluations. Thus, we are motivated by the potential for the COBALT algorithm to achieve the same level of accuracy with significantly fewer function evaluations. 

For completeness, we provide a brief summary of the problem; interested readers are referred to \cite{paulson19} for a more detailed description. The system of interest is for diauxic growth of a batch monoculture of \textit{Escherichia coli} on a glucose and xylose mixed media, which can be modeled using dynamic flux balance analysis (DFBA) \cite{hoffner13}. DFBA models are formulated as ordinary differential equations with embedded optimization problems that represent the metabolic network of the microorganisms, which means they are fundamentally a dynamic simulation with discrete events (i.e., a hybrid system). A general representation of DFBA models is
\begin{align} \label{eq:dfba-states}
\dot{\mathbf{s}}(t,\theta) = \mathbf{f}(t, \mathbf{s}(t, \theta), \mathbf{v}(\mathbf{s}(t, \theta), \theta), \theta), ~~ \mathbf{s}(t_0, \theta) = \mathbf{s}_0(\theta), 
\end{align}
where $\mathbf{v}(\mathbf{s}, \theta)$ is an element of the solution set of the metabolic flux balance model
\begin{align} \label{eq:flux-balance}
\mathbf{v}(\mathbf{s}, \theta) \in \argmax_\mathbf{v} &~~ q(\mathbf{v}, \mathbf{s}, \theta), \\\notag
\text{s.t.} &~~ \mathbf{A}(\theta)\mathbf{v} = \mathbf{0}, \\\notag
&~~ \mathbf{v}^{LB}(\mathbf{s}, \theta) \leq \mathbf{v} \leq \mathbf{v}^{UB}(\mathbf{s}, \theta),
\end{align}
$\mathbf{s}(t, \theta)$ denotes the state variables describing the extracellular environment at time $t$ for fixed parameters $\theta$ (e.g., concentration of biomass, substrates, and products) with initial conditions $\mathbf{s}_0(\theta)$; $\mathbf{v}$ denotes the metabolic fluxes that include both intracellular fluxes and exchange rates; $\mathbf{A}(\theta)$ is the  stoichiometric matrix of the metabolic network; $\mathbf{v}^{LB}(\mathbf{s}, \theta)$ and $\mathbf{v}^{UB}(\mathbf{s}, \theta)$ are the upper and lower bounds on the fluxes that depend on the extracellular concentration, respectively; $\mathbf{f}$ is a vector function that defines the rate of change of each component of $\mathbf{s}$ (specified by mass balances in the extracellular medium); and $q$ is the scalar function that represents the cellular objective function to be maximized. For the considered case of batch operation of an \textit{E. coli} fermentation reactor, we can specify \eqref{eq:dfba-states} in terms of three states, i.e., the concentrations of biomass, glucose, and xlyose. The flux balance model \eqref{eq:flux-balance} is constructed from wild-type \textit{E. coli} using the iJR904 metabolic network reconstruction \cite{reed03}, which contains 1075 reactions and 761 metabolites. The cellular objective was chosen to maximize the growth of the cells at every time point. We simulated this model using the DFBAlab toolbox \cite{gomez14} that utilizes Matlab's \texttt{ode15s} to integrate \eqref{eq:dfba-states} and CPLEX to solve the LP representation of \eqref{eq:flux-balance}. 

In this particular model, we have six unknown parameters $\theta$ that correspond to the maximum substrate uptake rates, saturation constants, and inhibition constants that appear $\mathbf{v}^{LB}(\mathbf{s}, \theta)$ and $\mathbf{v}^{UB}(\mathbf{s}, \theta)$. Noisy concentration measurements are available at eight time points $t \in \{ 5.5, 6.0, 6.5, 7.0, 7.25, 8.0, 8.25, 8.5 \}$:
\begin{align} \label{eq:dfba-meas}
\mathbf{y}_{i,j} = [\mathbf{s}(t_i, \theta)]_j + \mathbf{e}_{i, j}, ~~ \forall (i,j) \in \{1,\ldots, 8 \} \times \{ 1,2, 3\},
\end{align}
where $[\mathbf{s}]_j$ is the $j^\text{th}$ component of the state vector $\mathbf{s}$ that is composed of biomass ($j=1$), glucose ($j=2$), and xylose ($j=3$) and $\mathbf{e}_{i,j}$ is the measurement error for the $j^\text{th}$ state at the $i^\text{th}$ time point. Let $\hat{\mathbf{y}}_{i,j}(\theta) = [\mathbf{s}(t_i, \theta)]_j$ denote the model prediction of the $j^\text{th}$ state at the $i^\text{th}$ time point. Since the noise in the concentration sensors often depends on the absolute value of the measurement, we model $\mathbf{e}_{i,j}$ as zero-mean Gaussian random variables with a state-dependent standard deviation that is 5\% of the measured signal, i.e., 
\begin{align} \label{eq:dfba-noise}
\mathbf{e}_{i, j} \sim \mathcal{N}(0, \sigma_{i,j}^2(\theta)), ~~ \sigma_{i,j}(\theta) = 0.05| \hat{\mathbf{y}}_{i,j}(\theta) |, ~~ \forall (i,j) \in \{1,\ldots, 8 \} \times \{ 1,2,3\}.
\end{align}
The likelihood function is then specified by the collection of data and noise models in \eqref{eq:dfba-meas} and \eqref{eq:dfba-noise} as follows
\begin{align} \label{eq:dfba-likelihood}
p(\mathbf{Y}_\text{meas} | \theta) = \prod_{i=1}^8 \prod_{j=1}^3 \frac{1}{\sqrt{2 \pi \sigma^2_{i,j}(\theta)}} \exp\left( -\frac{(\mathbf{y}^\text{meas}_{i,j} - \hat{\mathbf{y}}_{i,j}(\theta))^2}{2\sigma^2_{i,j}(\theta)} \right),
\end{align}
where $\mathbf{Y}_\text{meas} = (\mathbf{y}^\text{meas}_{1,1}, \mathbf{y}^\text{meas}_{1,2}, \ldots, \mathbf{y}^\text{meas}_{8,3}) \in \mathbb{R}^{24}$ is the complete set of measured data. Our goal is to find the best parameter estimate by minimizing the negative log likelihood, i.e., $\theta_\text{MLE} = \argmin_\theta \{ -\log(p(\mathbf{Y}_\text{meas} | \theta) ) \}$. We can convert this problem to the standard form of \eqref{eq:grey-box-problem} by replacing $\theta \leftarrow x$ in the objective, rearranging, and removing constant terms that do not change the location of the minimum. The final description of the problem is provided in Appendix \ref{appendix:A}; notice the composite structure of $f$, which is highly nonlinear function of the outputs of the DFBA simulator that can be exploited by COBALT but is neglected by other purely black-box methods. 

\begin{figure}[tb!]
\centering
\includegraphics[width=0.6\linewidth]{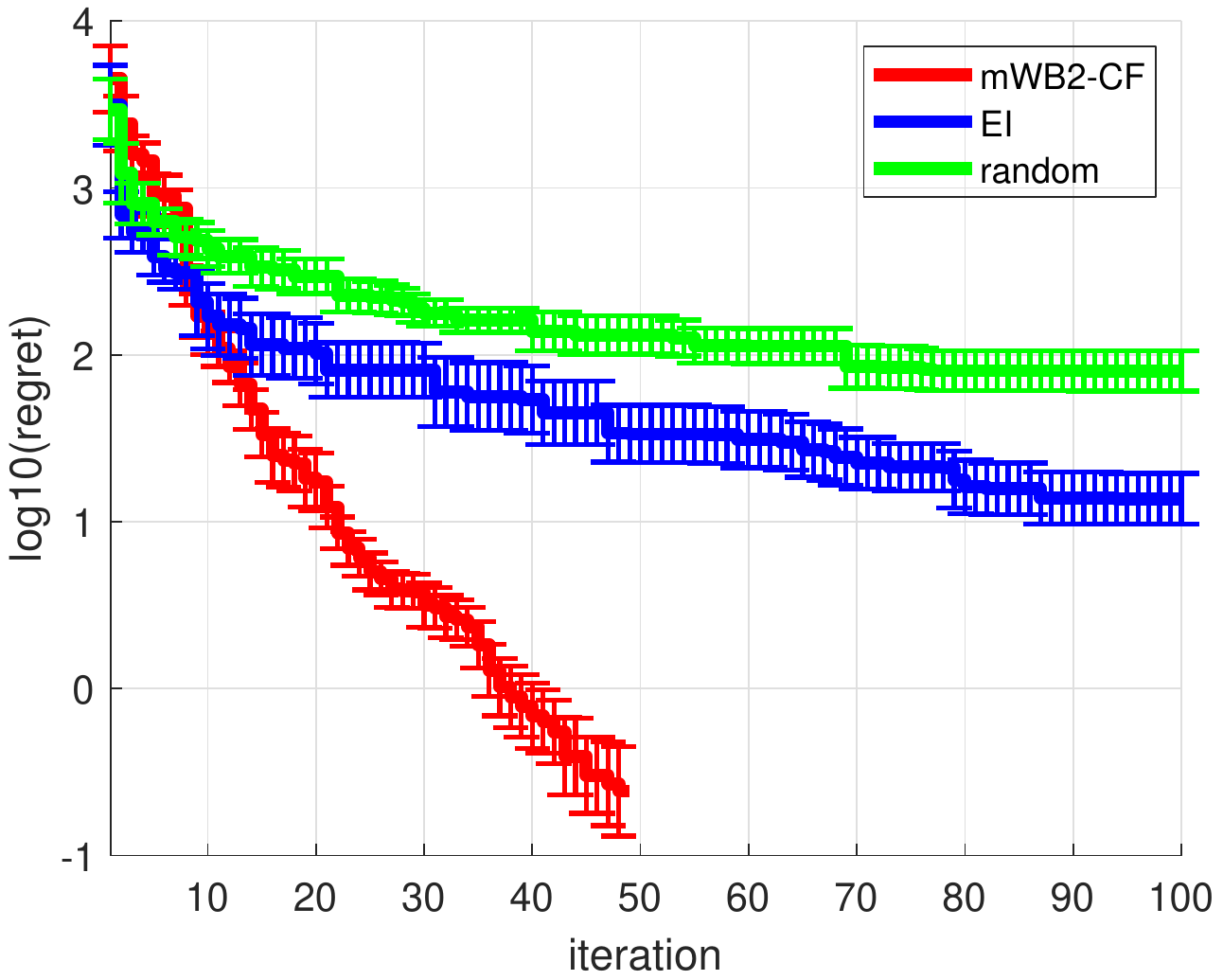}
\caption{Expected $\log10(\text{regret})$ for the DFBA-MLE problem, with approximate confidence region shown via error bars, estimated from 20 independent realizations.}
\label{fig:dfba}
\end{figure}

The performance of COBALT, relative to random search and EI-based BO, is shown in Fig. \ref{fig:dfba}. Note that, due to the complexity of these simulations, we reduced the number of replications to 20 and compared to a limited number of alternative methods. We can see that, although EI outperforms random search, it still reaches a relatively poor solution on average even when the number of iterations is increased to 100 objective evaluations. COBALT, however, is able to consistently converge to a low regret value (near the global minimum) in only 50 iterations. This is significantly fewer evaluations than the 1200 needed to construct the global surrogate in \cite{paulson19}. From this, we find that our results on this problem highlight two important points: (i) we can make significant improvements over purely black-box methods by exploiting structure whenever possible and (ii) we often needed fewer samples to find the most likely optimum point than to build a globally accurate surrogate model, which is especially important for computationally expensive simulators/models. 

%%%%%%%%%%%%%%%%%%%%%%%%%
\section{Conclusions and Future Work}
\label{sec:conclusions}

In this work, we present a novel algorithm for efficiently searching for the global optimum of constrained grey-box optimization problems. In particular, we extend the traditional black-box Bayesian optimization (BO) framework to handle composite functions of the form $f(x) = g(h(x))$ where $g(\cdot)$ is a known and differentiable white-box function and $h(\cdot)$ is an unknown vector-valued black-box function. These types of grey-box constraints, that involve a mixture of white- and black-box components, occur in a diverse set of engineering problems such as the calibration of high-fidelity simulators to experimental data, material and drug design, and maximum \textit{a posteriori} estimation of machine learning models with expensive-to-evaluate likelihoods. The proposed algorithm, which we refer to as COBALT, combines multivariate Gaussian process (GP) regression models (which are non-parametric and probabilistic in nature) with a novel expected utility (or acquisition) function that is subject to a chance constraint-based representation of the feasible region. Due to the general composite structure of the grey-box objective and constraint functions, we cannot simply maximize the constrained acquisition function using standard optimization methods, as is the case in traditional BO. Instead, we propose a sample average approximation of the intractable expectation-based objective to convert it into a deterministic expression whose first- and second-order derivatives can be efficiently computed. In addition, we develop a simple moment-based approximation of the chance constraints, so that the overall maximization problem can be efficiently solved with state-of-the-art nonlinear programming solvers. 

To demonstrate the effectiveness of the proposed algorithm, we compare COBALT to traditional BO on a set of seven diverse test optimization problems. We observe that COBALT outperforms BO, with respect to the average regret (i.e., difference between the estimated optimum at each iteration and the true global minimum), in all considered cases. The final test problem is a realistic engineering problem that is focused on estimating parameters in a genome-scale bioreactor model using experimental measurements. This problem, which was recently tackled for the first time in \cite{paulson19} using a custom surrogate-based optimization paradigm, required around 1200 expensive model simulations to build a surrogate model that was accurate enough to (approximately) solve the original parameter estimation problem. When solving this problem using COBALT in this work, we found that accurate solutions (i.e., within $<1\%$ of the best known optimum) could be found in 50 iterations on average, implying a factor of 24 less total number of expensive model simulations. Furthermore, COBALT found solutions that were two orders of magnitude better than standard BO found in double the number of iterations (100 total). We believe that the results presented in this paper indicate that a promising path toward accelerating the convergence of black-box optimization methods is exploiting known problem structure as much as possible. In addition, we find that the results highlight the fundamental difference between constructing a globally predictive surrogate model and locating the minimum of a constrained function. It turns out the latter problem is easier than the former, as one can simply ignore parts of the design space that are not likely to produce ``good'' solutions. 

Although promising results have been obtained in Section \ref{sec:numerical-example}, there remain several important directions for future work that could improve the numerics and theoretical understanding of COBALT. We briefly discuss three possible directions below. 

\subsection{Extensions to other acquisition functions}

Throughout this work, we focus on modified forms of the expected improvement (EI) acquisition function, which was mostly motivated by results in \cite{astudillo19}. However, it is well-known that EI works best under the assumption that each observation of the objective and/or constraints has low-noise relative to the true value of the function. This assumption is not necessarily satisfied in several application domains including drug discovery, medical diagnostics, biosurveillance, and molecular simulation models. Two particularly interesting alternatives to EI in such cases are the knowledge gradient (KG) \cite{frazier09} and predictive entropy search (PES) \cite{hernandez14} acquisition functions. Instead of relying on a selected incumbent value to guide the selection of the next sample, KG and PES more heavily rely on the posterior distribution of the surrogate model itself and thus are less sensitive to noise. If more exploration is desired, then one could alternatively select the lower confidence bound (LCB) \cite{srinivas09} acquisition function or apply Thompson sampling (TS) \cite{bradford18}. To the best of our knowledge, KG, PES, LCB, or TS have not been extended to the case of composite functions. Since KG and PES are computationally intensive in the case of standard BO, we expect new tractable approximation methods will need to be developed to apply under our proposed constrained grey-box setting \eqref{eq:grey-box-problem}.  

\subsection{Improving the enrichment sub-problem optimization scheme}
\label{subsec:improving-opt-subproblem}

As shown in \eqref{eq:saa}, we had to invoke two main approximations (i.e., sample average approximation of the objective and moment-based approximation of the chance constraints) to make the enrichment sub-problem \eqref{eq:optimize-acqusition} tractable. Although a large amount of analysis has been performed on these approximations (see, e.g., Remark \ref{rem:2}), we have not explicitly made use of this theory in COBALT. In particular, we do not adaptively set the number of random samples $M$ to ensure the solution of the approximate problem is within a specified tolerance of the true solution. This is mainly due to the added complexity of solving many SAA problems at each iteration of COBALT. One potential strategy for reducing the cost of these SAA problems is to take advantage of parallel or distributed optimization methods, such as the alternating direction method of multipliers (ADMM) \cite{boyd11}, that can achieve significant computational speedup and/or reduced communication between sub-problems. In addition, the stochastic gradient ascent method presented in \eqref{eq:sga} can be modified to ensure robust convergence and even outperform SAA for certain classes of problems \cite{nemirovski09} (for convex constraint sets $\mathcal{X}_n$). As such, developing new ways to achieve similar convergence guarantees in the presence of non-convex sets $\mathcal{X}_n$ is also an interesting direction for future work. Finally, we note that the EI-CF and mWB2-CF acquisition functions both have the potential to have many local solutions that may provide useful information at each iteration. Whenever multiple function evaluations can be run in parallel at each iteration, we can use these samples to develop a batch version of COBALT (see, e.g., \cite{liu21} for more details on batch BO). 

\subsection{Handling high-dimensional problems using trust regions}

As discussed in detail in the introduction, COBALT is focused on \textit{global} optimization of constrained grey-box models, which requires a global surrogate model to constructed at every iteration. As the number of inputs to the black-box portion of the model grows, the construction of this surrogate gets significantly more challenging, implying COBALT (as presented in this paper) may not perform well on particularly large-scale problems. The most straightforward way to overcome this challenge is to apply dimensionality reduction techniques before constructing the GP model \cite{bouhlel16}; building specialized methods that can be done simultaneously with the GP repression step may be an interesting to pursue. An alternative direction would involve combining COBALT with established trust region methods, e.g., \cite{eason16, eason18, bajaj18}, which only construct surrogate models over a local region of the design space. Although there are established convergence results for various trust region methods, they can only guarantee convergence to a local optimum and rely on several tuning parameters that may be difficult to select before running the algorithm. In addition, most trust region methods are based on deterministic models and may require fully linear models \cite{wild08} to be constructed using, e.g., linear/quadratic interpolation, which means the number of required evaluations at each iteration scales with the size of the input dimension. A trust region BO algorithm, TuRBO \cite{eriksson19}, was recently developed that addresses the local versus global tradeoff by running several independent local models and uses an implicit multi-arm bandit approach to decide which local model should be allocated samples. Ideas from TuRBO could inform a local trust region-based strategy for COBALT, which may prove to be beneficial in the context of large-scale problems. 
% switching algorithm that seems interesting: https://arxiv.org/pdf/1805.08610.pdf

%Bound regret
%% https://proceedings.neurips.cc/paper/2018/file/41f860e3b7f548abc1f8b812059137bf-Paper.pdf
%
%EI-CF paper has some theory
%% discuss that

% Other acquisition functions such as knowledge gradient. Easy to extend to LCB, which we will explore more in future. TS is simple. 
%Sparse GP. 
% Ways to use SGD still in the subproblem solution or ways to handle chance constraints more exactly. Distributed optimization also. 
% Trust regions to scale better. 
%Noisy function evaluations using input-dependent variance. We did not focus much, but we can handle noisy simulations which are present in many situations. However, this is only handled through $\sigma_v^2$, which is assumed to be constant -- want to extend. 
% Prove convergence / properties in the constrained case. [there is some initial theory for EI-CF]. 

%%%%%%%%%%%%%%%%%%%%%%%%%
\appendix
\numberwithin{equation}{section}
\section{Appendix: Representation of test problems in standard form}
\label{appendix:A}

\subsection{Goldstein-Price}

We propose a modified formulation of the Goldstein-Price function that can be formulated as a grey-box optimization problem \eqref{eq:grey-box-problem} as follows
\begin{align} \label{eq:GSP}
\min_{x, y, z} &~~ (1 + (x_1 + x_2 +1)^2(19 - 14 x_1 + 3x_1^2 + y_1)) \\\notag
&~~~ \cdot (30 + y_2(18-32x_1 + 12x_1^2 + 48 x_2 - 36 x_1 x_2 + 27x_2^2)),  \\\notag
\text{s.t.} &~~ y_1 = d_1(z) := -14 z_2 + 6 z_1 z_2 + 3 z_2^2, \\\notag
&~~ y_2 = d_2(z) := 2 z_1^2 - 3 z_2^2, \\\notag
&~~ z = x, \\\notag
&~~ -2 \leq x_i \leq 2, ~~~~~ \forall i \in \{ 1, 2 \}.
\end{align}
The global minimum is equal to $3$ with $x^\star = [0, -1]^\top$.

\subsection{Rastrigin}

We propose a modified formulation of the Rastrigin function that can be formulated as a grey-box optimization problem \eqref{eq:grey-box-problem} as follows
\begin{align} \label{eq:Rastrigin}
\min_{x, y, z} &~~ 30 + x_1^2 - 10\cos(2\pi x_1) + x_2^2 - 10\cos(2\pi x_2) + y_1, \\\notag
\text{s.t.} &~~ y_1 = d_1(z) := z_1^2 - 10\cos( 2 \pi z_1), \\\notag
&~~ z_1 = x_3, \\\notag
&~~ -5.12 \leq x_i \leq 5.12, ~~~~~ \forall i \in \{ 1, 2,3 \}.
\end{align}
The global minimum is equal to $0$ with $x^\star = [0, 0, 0]^\top$. 

\subsection{Rosenbrock}

We propose a modified formulation of the Rosenbrock function that can be formulated as a grey-box optimization problem \eqref{eq:grey-box-problem} as follows
\begin{align} \label{eq:Rosenbrock}
\min_{x, y, z} &~~ \sum_{i=1}^3 (100 y_i^2  + (1-x_i)^2) + 100(x_5 - x_4^2) + y_4 \\\notag
&~~~~~~~~ + 100(x_6 - x_5^2) + (1 - x_5)^2, \\\notag
\text{s.t.} &~~ y_i = d_i(z) := z_{i+1}^2 - z_i^2, ~~~~~ \forall i \in \{ 1, 2, 3 \} , \\\notag
&~~ y_4 = d_4(z) := (1 - z_4)^2, \\\notag
&~~ z=x, \\\notag
&~~ -2 \leq x_i \leq 2, ~~~~~ \forall i \in \{ 1, 2, 3, 4 \}.
\end{align}
The global minimum is $0$ with $x^\star = [0, 0, 0, 0, 0, 0]^\top$.

\subsection{Toy-Hydrology}

We propose a modified formulation of the Toy-Hydrology function that can be formulated as a constrained grey-box optimization problem \eqref{eq:grey-box-problem} as follows
\begin{align} \label{eq:ToyHydrology}
\min_{x, y, z} &~~ x_1 + x_2, \\\notag
\text{s.t.} &~~ g_1(x, y) := 1.5 - x_1 - 2x_2 - 0.5\sin(-4\pi x_2 + y_1) \leq 0, \\\notag
&~~ g_2(x, y) := x_1^2 + x_2^2 - 1.5 \leq 0, \\\notag
&~~ y_1 = d_1(z) := 2 \pi z_1^2, \\\notag
&~~ z_1 = x_1, \\\notag
&~~ 0 \leq x_i \leq 1, ~~~~~ \forall i \in \{ 1, 2 \}, 
\end{align}
The global minimum is $0.5998$ with $x^\star = [0.1951, 0..4047]^\top$.

\subsection{Rosen-Suzuki}

We propose a modified formulation of the Rosen-Suzuki function that can be formulated as a constrained grey-box optimization problem \eqref{eq:grey-box-problem} as follows
\begin{align} \label{eq:RosenSuzuki}
\min_{x, y, z} &~~ x_1^2 + x_2^2 + x_4^2 - 5x_1 - 5x_2 + y_1, \\\notag
\text{s.t.} &~~ g_1(x, y) := -(8 - x_1^2 - x_2^2 - x_3^2 - x_4^2 - x_1 + x_2 - x_3 + x_4) \leq 0, \\\notag
&~~ g_2(x, y) := -(10 - x_1^2 - 2x_2^2 - y_2 + x_1 + x_4) \leq 0, \\\notag
&~~ g_3(x, y) := -(5 - 2x_1^2 - x_2^2 - x_3^2 - 2x_1 + x_2 + x_4) \leq 0, \\\notag
&~~ y_1 = d_1(z) := 2z_1^2 - 21z_1 + 7z_2, \\\notag
&~~ y_2 = d_2(z) := z_1^2 + 2z_2^2, \\\notag
&~~ z_1 = x_3, \\\notag
&~~ z_2 = x_4, \\\notag
&~~ -2 \leq x_i \leq 2, ~~~~~ \forall i \in \{ 1,\ldots, 4 \}.
\end{align}
The global minimum is $-44$ with $x^\star = [0, 1, 2, -1]^\top$.

\subsection{Colville} \label{subappendix:colville}

We propose a modified formulation of the Colville function that can be formulated as a constrained grey-box optimization problem \eqref{eq:grey-box-problem} as follows
\begin{align} \label{eq:Colville}
\min_{x, y, z} &~~ 5.3578 x_3^2 + y_1, \\\notag
\text{s.t.} &~~ g_1(x, y) := y_2 - 0.0000734x_1x_4 - 1 \leq 0, \\\notag
&~~ g_2(x, y) := 0.000853007x_2x_5 + 0.00009395x_1x_4 - 0.00033085x_3x_5 - 1 \leq 0, \\\notag
&~~ g_3(x, y) := y_4 - 0.30586 (x_2x_5)^{-1} x_3^2 - 1 \leq 0, \\\notag
&~~ g_4(x, y) := 0.00024186x_2x_5 + 0.00010159x_1x_2 + 0.00007379x_3^2 - 1 \leq 0, \\\notag
&~~ g_5(x, y) := y_3 - 0.40584(x_5)^{-1}x_4 - 1 \leq 0, \\\notag
&~~ g_6(x, y) := 0.00029955x_3x_5 + 0.00007992x_1x_3 + 0.00012157x_3x_4 - 1 \leq 0, \\\notag
&~~ y_1 = d_1(z) := 0.8357z_1z_4 + 37.2392z_1, \\\notag
&~~ y_2 = d_2(z) := 0.00002584z_3z_4 - 0.00006663z_2z_4, \\\notag
&~~ y_3 = d_3(z) := 2275.1327(z_3z_4)^{-1}- 0.2668(z_4)^{-1}z_1, \\\notag
&~~ y_4 = d_4(z) := 1330.3294(z_2z_4)^{-1} - 0.42(z_4)^{-1}z_1, \\\notag
&~~ z_1 = x_1, \\\notag
&~~ z_2 = x_2, \\\notag
&~~ z_3 = x_3, \\\notag
&~~ z_4 = x_5, \\\notag
&~~ 78 \leq x_1 \leq 102, \\\notag
&~~ 33 \leq x_2 \leq 45, \\\notag
&~~ 27 \leq x_i \leq 45, ~~~~~ \forall i \in \{ 3,4,5 \}.
\end{align}
The global minimum is $10122.7$ with $x^\star = [78, 33, 29.998, 45, 36.7673]^\top$.

\subsection{DFBA-MLE}

The DFBA-MLE problem seeks to minimize the negative log of the likelihood function that was derived in \eqref{eq:dfba-likelihood}. After some algebraic manipulations, we arrive at the following grey-box optimization problem in the form of \eqref{eq:grey-box-problem}
\begin{align} \label{eq:dfba-substitutions}
\min_{x, y, z} &~~ \sum_{i=1}^8 \sum_{j=1}^3 \log(0.0025 y_{i+3(j-1)}^2) + \frac{( \mathbf{y}^\text{meas}_{i,j} - y_{i+3(j-1)}  )^2}{0.0025 y_{i+3(j-1)}^2}, \\\notag
\text{s.t.} &~~ y_{i+3(j-1)} = d_{i+3(j-1)}(z) := \hat{\mathbf{y}}_{i,j}(z), ~~~~~ \forall (i,j) \in \{1,\ldots, 8 \} \times \{ 1,2,3\}, \\\notag
&~~ z = x, \\\notag
&~~ 9.45 \leq x_1 \leq 11.55, \\\notag
&~~ 0.0024 \leq x_2 \leq 0.0030, \\\notag
&~~ 5.4 \leq x_3 \leq 6.6, \\\notag
&~~ 0.0149 \leq x_4 \leq 0.0182, \\\notag
&~~ 0.0045 \leq x_5 \leq 0.0055, \\\notag
&~~ 12.2727 \leq x_6 \leq 15,
\end{align}
where we have substituted $z = x = \theta$. The exact global minimum for this problem is unknown and so was estimated from the minimum objective value obtained across all runs of every algorithm. To check that this identified solution is likely to be near the global, we verified that the resulting parameter estimate gave high-quality predictions that result in a very large likelihood value.  

\section{Appendix: Overview of Computational Statistics for COBALT}
\label{appendix:B}

In this section, we provide a detailed discussion of the computational cost of the proposed COBALT method relative to standard constrained Bayesian optimization (CBO) methods. The time it takes to execute iteration $n$ of COBALT, excluding the evaluation cost in Step 4 of Algorithm \ref{alg:gp-sampling}, is given by
\begin{align} \label{eq:time-cobalt}
t_{\text{COBALT}, n} = \textstyle \sum_{i=1}^{n_y} t_{\text{GP}, n, i} + t_{\text{GreyOpt}, n},
\end{align}
where $t_{\text{GP}, n, i}$ is the time it takes to train the GP model for component $i \in \{1,\ldots, n_y\}$ of $d(z)$ and $t_{\text{GreyOpt}, n}$ is the time required to solve the grey-box enrichment sub-problem \eqref{eq:optimize-acqusition}. 
On the other hand, the time it takes to execute iteration $n$ of a standard CBO algorithm (again excluding function evaluation cost) is equal to
\begin{align} \label{eq:time-cbo}
t_{\text{CBO}, n} = \textstyle \sum_{i=1}^{U + 1} t_{\text{GP}, n, i} + t_{\text{BlackOpt}, n},
\end{align}
where $U+1$ is the total number of black-box functions in the formulation \eqref{eq:simplified-grey-box} (composed of $U$ unknown constraints and 1 unknown objective that is equivalent to the grey-box formulation \eqref{eq:grey-box-problem}) and $t_{\text{BlackOpt}, n}$ is the time required to solve the black-box enrichment sub-problem. 

Let $t_{\text{Eval}, n}$ denote the time required to evaluate the unknown function at iteration $n$. In the case that $t_{\text{Eval}, n}$ is much larger than $t_{\text{COBALT}, n}$ and $t_{\text{CBO}, n}$, then the total time to execute $N$ steps is $t_\text{Total} \approx \sum_{n=1}^N t_{\text{Eval}, n}$, which is the same for COBALT and CBO. However, it is important to look at the absolute costs of training the GP and optimizing the acquisition function to determine how large $t_{\text{Eval}, n}$ must become for this assumption to be valid. 

It is interesting to look at the difference between \eqref{eq:time-cobalt} and \eqref{eq:time-cbo} assuming that the same GP training tool is used
\begin{align} \label{eq:time-difference}
t_{\text{COBALT}, n} - t_{\text{CBO}, n} = \textstyle \sum_{i=1}^{n_y-U-1} t_{\text{GP}, n, i} + t_{\text{GreyOpt}, n} - t_{\text{BlackOpt}, n}.
\end{align}
An immediate observation from this expression is that COBALT is not necessarily more expensive than CBO. In particular, the first term will be negative when $n_y < U + 1$, which can happen in problems that involve a significant number of constraints (not uncommon in engineering systems). In addition to the GP training time, the time required to numerically optimize the sub-problem plays an important role. Thus, it is difficult to perform a comprehensive and fair comparison since the values of $\{ t_{\text{GP}, n, i}, t_{\text{GreyOpt}, n}, t_{\text{BlackOpt}, n} \}$ depend strongly on the implementation details (including the choice of programming language and numerical optimization procedures). We discuss some of these complexities in the following three sections to provide some insight into the choices made in our initial implementation. Note that there remains room for improvement in terms of improving the efficiency of the selected algorithms; COBALT can directly benefit from any new developed algorithms that accelerate GP training or the numerical solution of the sub-problems.  

%A  is thus difficult to perform since there is not a single value for , which depend 
%
%to the time it takes to train the GP models, the time required to 
%
%whenever $U > 1$ -- in fact
%
%It is important to note that there is not a single value for , as they depend strongly on the details of the implementation (including the choice of programming language and optimization methods for solving the MLE problem and the acquisition sub-problems). 

\subsection{Cost of Gaussian process regression}

\begin{figure}[hbt]
\centering
\includegraphics[width=0.6\linewidth]{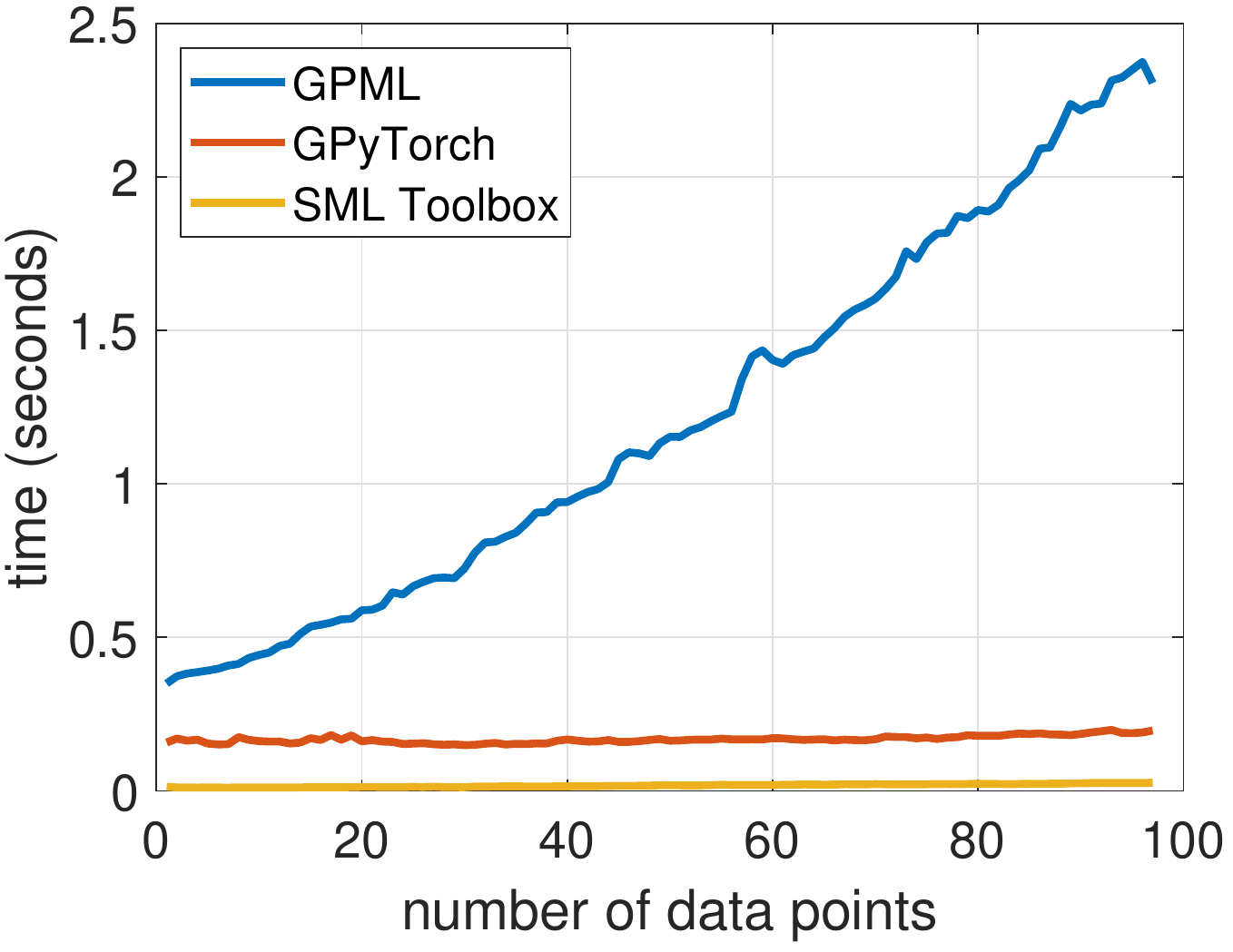}
\caption{CPU time required to train a Gaussian process model versus the number of measured data points for three different packages.}
\label{fig:gp-fit-time}
\end{figure}

Gaussian process (GP) regression has been an active area of research for more than two decades, with several packages available for performing the steps outlined in Section \ref{subsec:GPsingle}. To get an estimate of the training time, we analyze the performance of three commonly used packages on a simple test problem. In particular, we compare (i) GPML \cite{rasmussen10}, which was the toolbox we used in COBALT, (ii) the \texttt{fitgpr} function implemented as a part of the Statistics and Machine Learning (SML) toolbox in Matlab, and (iii) the GPyTorch toolbox implemented in Python. Note that all three of these packages solve the MLE problem \eqref{eq:MLE} using different strategies; we used the default settings for simplicity, though additional tuning could be done to improve performance in each of these packages. The CPU time\footnote{All computational experiments carried out in Appendix \ref{appendix:B} were run on a laptop with 16 GB of RAM and a 1.8GHz Intel i7 processor.} versus the number of training data points for the three considered GP packages is shown in Fig. \ref{fig:gp-fit-time}. We see that GPML takes the longest time, which is partially due to significant overheard incurred by inverting the kernel matrix that has not been optimized. However, the overall process takes a maximum of 2.4 seconds for 100 training data points, which is a reasonable cost for an expensive function (that could take many hours or more to evaluate). 

\begin{figure}
     \centering
     \begin{subfigure}[b]{0.6\textwidth}
         \centering
         \includegraphics[width=\textwidth]{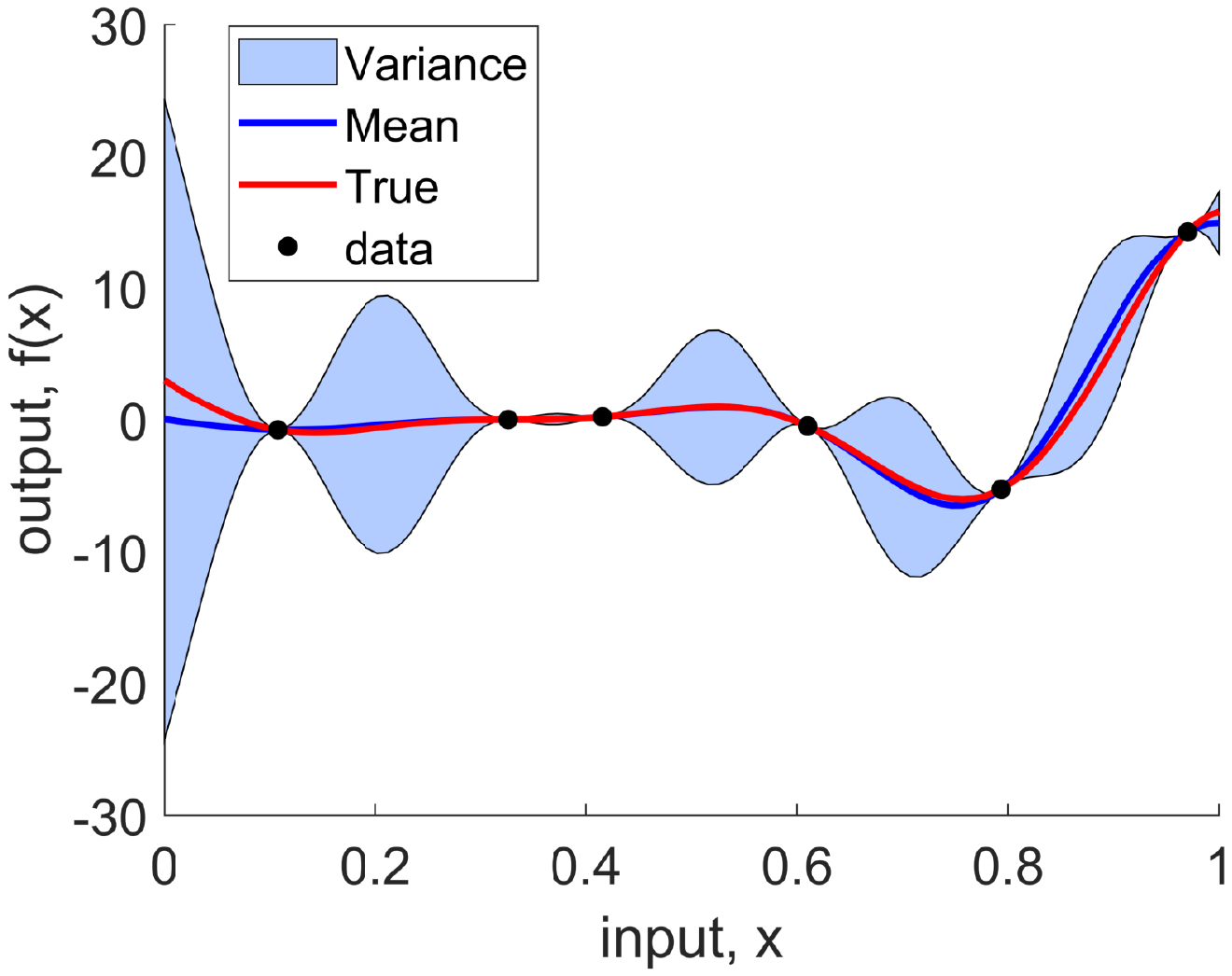}
         \caption{GPML}
     \end{subfigure}
     \hfill
     \begin{subfigure}[b]{0.6\textwidth}
         \centering
         \includegraphics[width=\textwidth]{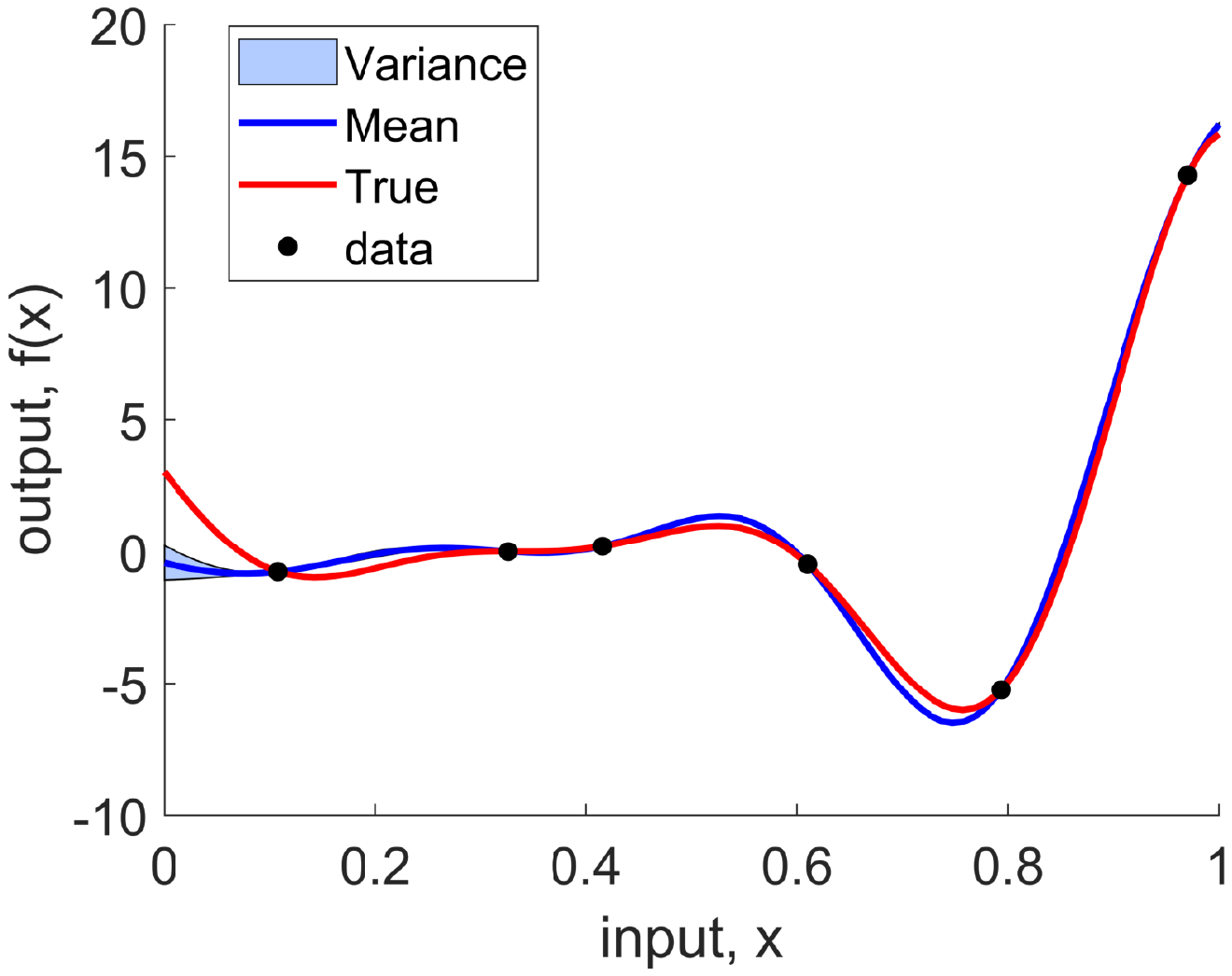}
         \caption{GPyTorch}
     \end{subfigure}
     \hfill
     \begin{subfigure}[b]{0.6\textwidth}
         \centering
         \includegraphics[width=\textwidth]{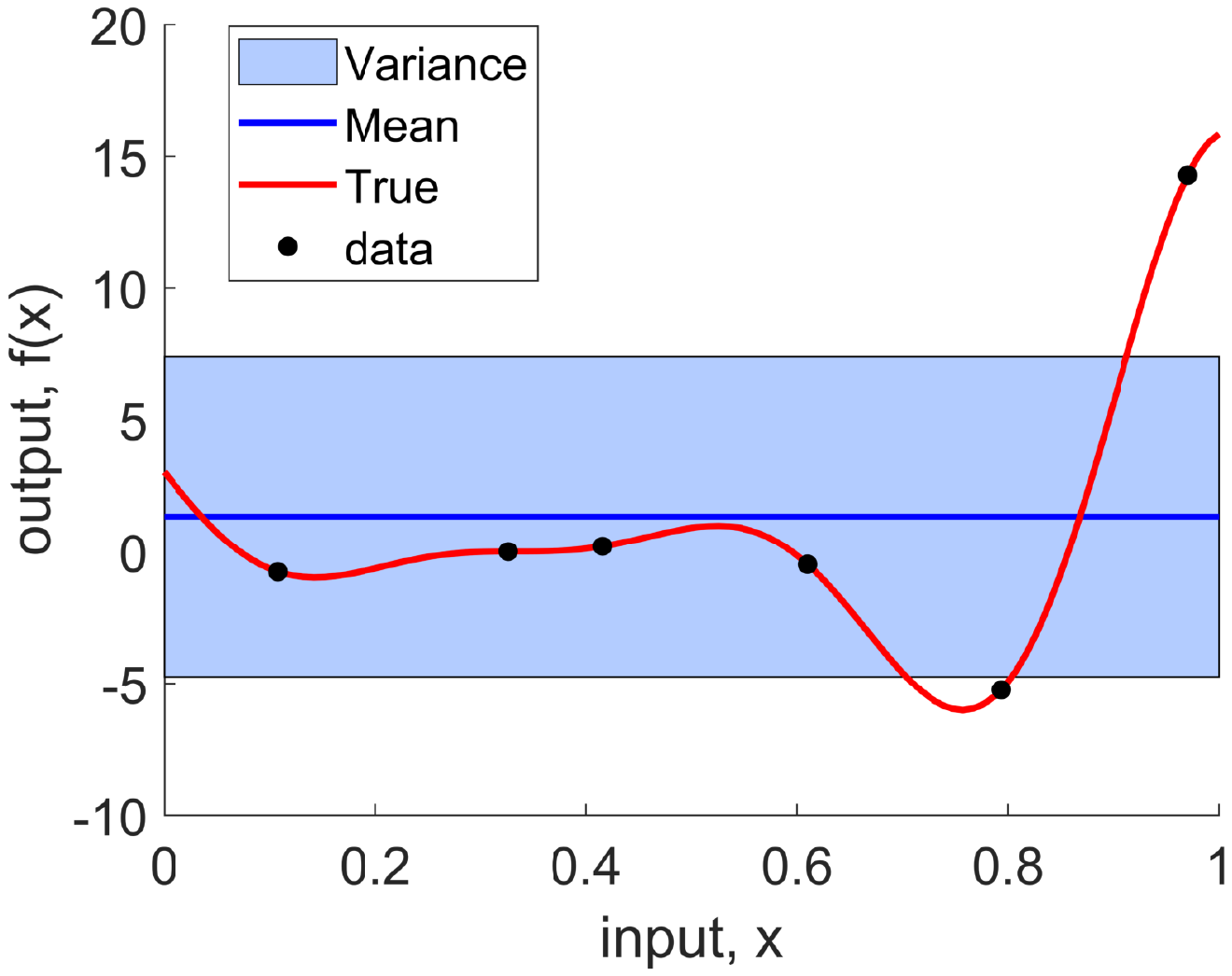}
         \caption{SML Toolbox}
     \end{subfigure}
        \caption{Comparison of the Gaussian process model predictions obtained with three different packages: (a) GPML, (b) GPyTorch, and (c) SML Toolbox. The solid blue line represents the mean prediction and blue cloud denotes the confidence region ($\pm1$ standard deviation). The true unknown function is shown with a solid red line and the available data is shown with black dots.}
        \label{fig:gp-prediction}
\end{figure}

It is interesting to note that GPyTorch and SML appear to scale more favorably with respect to the amount of training data. To investigate this further, we looked at the quality of the trained GP models compared to the true unknown function, which is shown in Fig. \ref{fig:gp-prediction}. From these plots, we see that GPML yields the best predictions, with the true function being fully contained within the estimated confidence region. GPyTorch and SML, on the other hand, under- and over-predict the true uncertainty in the function, respectively, which can be attributed to identifying high sub-optimal solutions to \eqref{eq:MLE}. As expected, the speed afforded by alternative GP regression tools may come at the cost of accuracy. Based on these results, we opted to use GPML in our implementation of COBALT since we found it to be quite reliable (and not overly expensive) in the low-data regime. 

\subsection{Cost of optimizing the proposed acquisition sub-problem}

\begin{figure}[hbt]
\begin{subfigure}{\linewidth}
\includegraphics[width=0.49\linewidth]{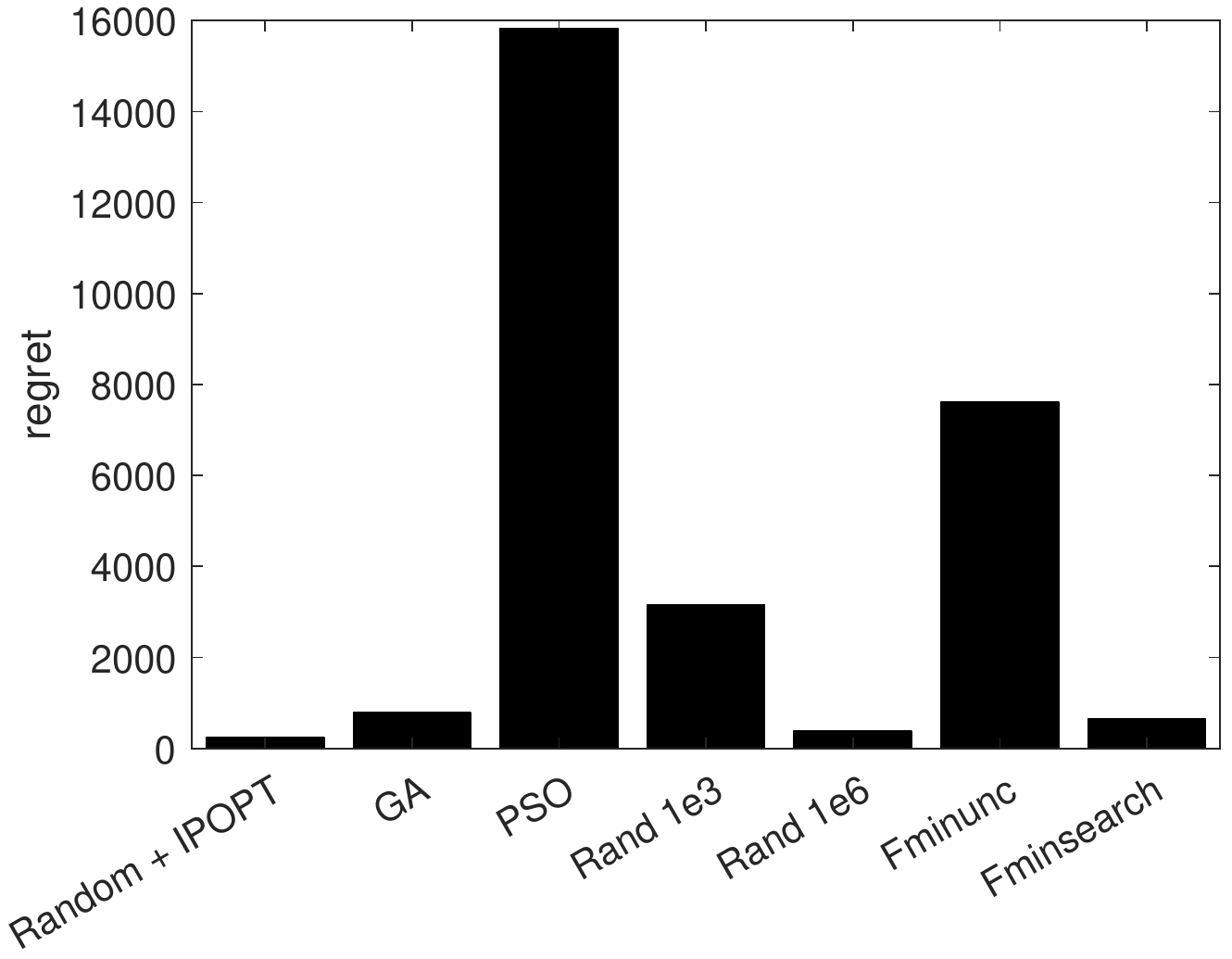}\hfill
\includegraphics[width=0.49\linewidth]{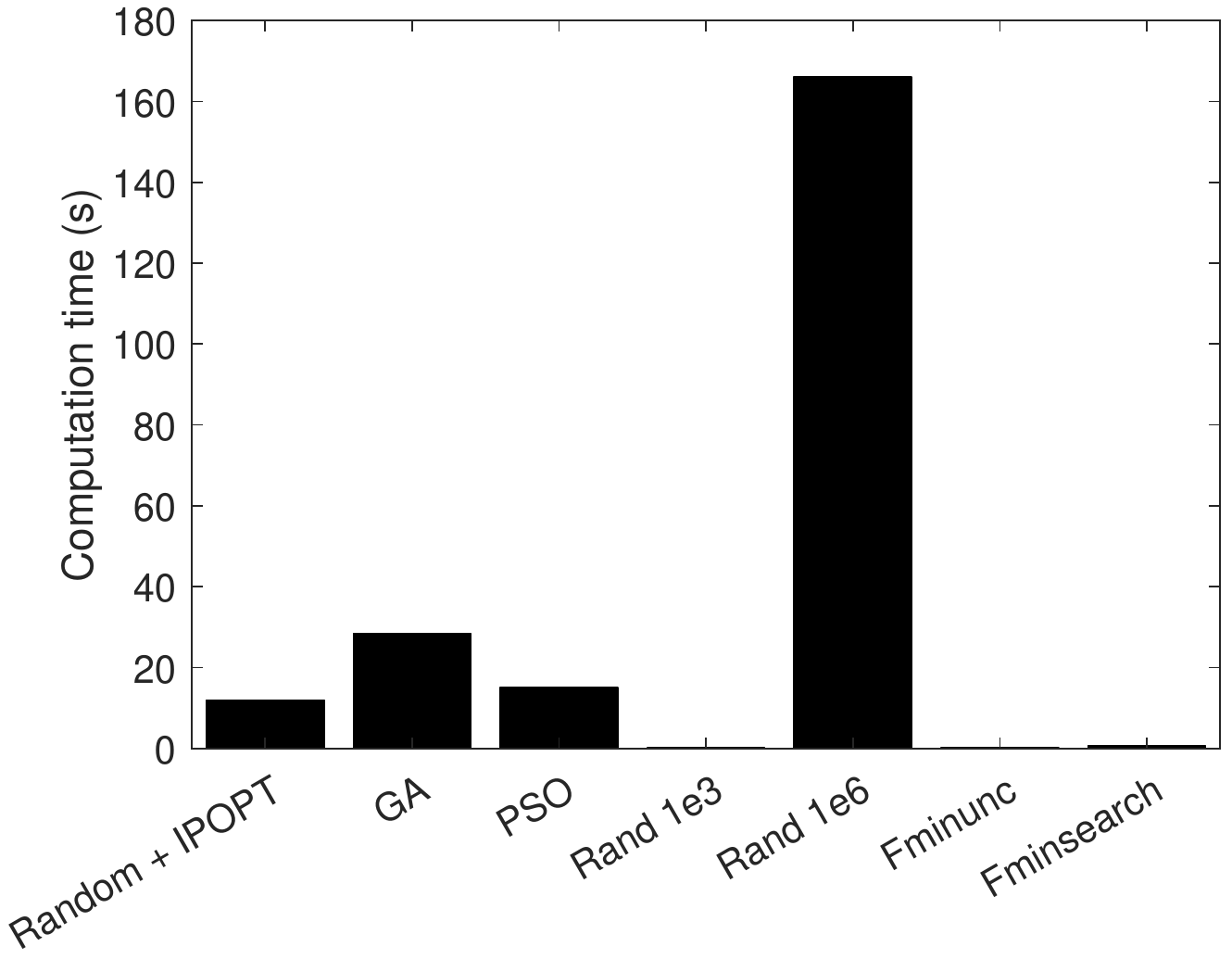}\hfill
\end{subfigure}\par\medskip
\caption{Comparison of the regret (left) and average CPU time (right) for solving an acquisition sub-problem (estimated by averaging over 20 independent replications) with seven different optimization algorithms.}
\label{fig:subopt}
\end{figure}

There are wide-variety of available optimization methods for solving the acquisition sub-problems, which are generally non-convex nonlinear programs. For reasons discussed in Section \ref{sec:cobalt}, we selected the local gradient-based solver IPOPT that we globalize with a random search approach (Random+IPOPT). To demonstrate the value of the proposed approach, we compared it to several alternative methods on the 6-dimensional Rosenbrock test problem (with 30 randomly drawn samples). We considered the following alternatives that are all readily available in Matlab: the genetic algorithm (GA) \texttt{ga}, particle swarm optimization (PSO) \texttt{pso}, random search (Rand), the quasi-Netwon solver \texttt{fminunc}, and the simplex-based search algorithm \texttt{fminsearch}. The regret, which is equal to $f(x_\text{estimated})-f(x^\star)$ where $f(x_\text{estimated})$ is the best solution returned by the optimizer and $f(x^\star)$ is the true global solution, and CPU time averaged over 50 replicates is shown in Fig. \ref{fig:subopt}. We see that our approach demonstrates a good tradeoff in that it involves a modest computational cost while providing consistently low regret values. However, as discussed in Section \ref{subsec:improving-opt-subproblem}, our proposed numerical solution method to the acquisition sub-problem is by no means optimized and there remains several interesting directions for future work. It should also be noted that further gains will be achieved in the numerical optimization of the sub-problem by accelerating the GP prediction (as this translates to cheaper evaluation of the acquisition and constraint functions). 

\subsection{Comparison of total CPU time per iteration for COBALT and CBO}

\begin{figure}[hbt]
\centering
\includegraphics[width=0.6\linewidth]{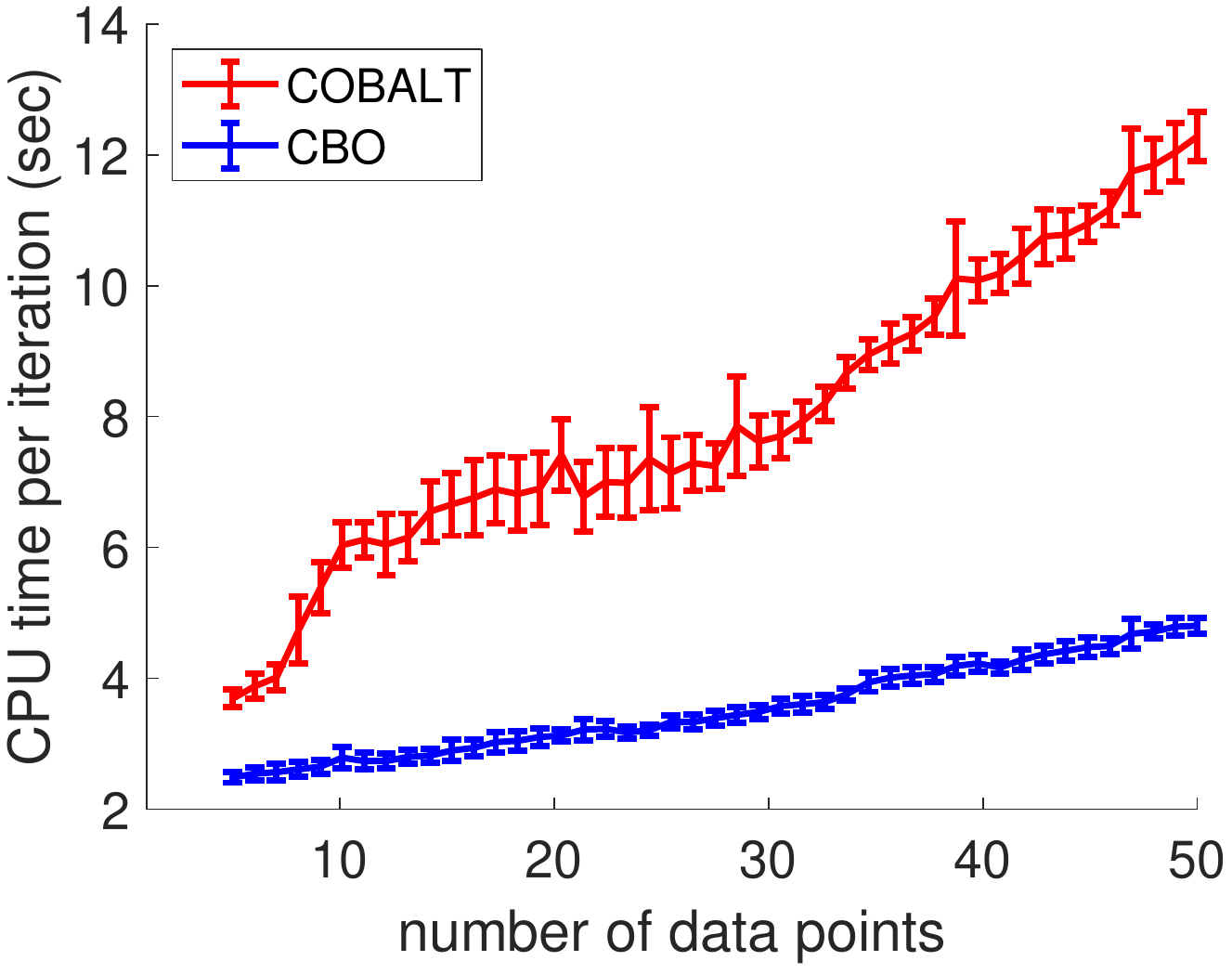}
\caption{Comparison of the total CPU time per iteration for COBALT and CBO on the Rosenbrock test problem under the same GP training and acquisition optimization procedures.}
\label{fig:time-cobalt-overall}
\end{figure}

As a final step, we look at the total execution time $t_{\text{COBALT}, n}$ and $t_{\text{CBO}, n}$ on the 6-dimensional Rosenbrock test problem. Note that we compare these two under the same GP training and acquisition optimization procedures discussed in the previous two sections -- it would not be fair to compare to \texttt{bayesopt}, which internally uses substantially different methods. 
The total CPU time per iteration versus the number of data points, averaged over 50 replications of the random initialization, is shown in Fig. \ref{fig:time-cobalt-overall}. Based on the analysis in \eqref{eq:time-difference}, we expect $t_{\text{COBALT}, n} > t_{\text{CBO}, n}$ since $n_y=4$ and $U=0$ for the Rosenbrock problem. The difference grows slightly as the number of data points $n$ increases since $t_{\text{GreyOpt}, n} > t_{\text{BlackOpt}, n}$ in this case due to the larger number of GP models; however, the growth is not substantial relative to the GP training cost, which dominates the cost at every iteration. For $n=50$, we see around an 8 second increase in cost for COBALT compared to CBO on this problem. We obtained similar results for the other test problems. 

Note that a fairly conservative estimate for $t_{\text{COBALT}, n}$ can easily be obtained using \eqref{eq:time-cobalt} with $t_{\text{GP}, n, i} \approx 5$ seconds and $t_{\text{GreyOpt}, n} \approx 30$ seconds. Therefore, we expect that COBALT would be the preferred option for any problem that has an internal component that takes on the order of minutes (or longer) and a budget of around $N= O(100)$. 

%% References
\bibliographystyle{elsarticle-num}           % Include this if you use bibtex
\bibliography{references}

\end{document}